\documentclass[11pt]{amsart}
\usepackage{amssymb,amsmath,amsthm}
\newcommand{\leg}[2]{\genfrac{(}{)}{}{}{#1}{#2}}

\newtheorem{theorem}{Theorem}
\newtheorem{lemma}[theorem]{Lemma}
\newtheorem{corollary}[theorem]{Corollary}

\theoremstyle{remark}

\newtheorem*{remark}{Remark}

\numberwithin{theorem}{section} \numberwithin{equation}{section}

\newcommand{\R}{\mathbb{R}}
\newcommand{\C}{\mathbb{C}}

\newcommand{\Z}{\mathbb{Z}}
\newcommand{\N}{\mathbb{N}}

\newcommand{\lcm}{{\text {\rm lcm}}}

\newcommand{\re}{\textnormal{Re}}

\def\H{\mathbb{H}}

\begin{document}
\title[On the explicit construction of  higher deformations of partition statistics]{On the explicit construction of  higher deformations of partition statistics} 

\author{Kathrin Bringmann}
\address{School of Mathematics\\University of Minnesota\\ Minneapolis, MN 55455 \\U.S.A.}
\email{bringman@math.umn.edu}      

\date{\today}
\maketitle
\begin{abstract}
The modularity of the partition generating function has many important consequences, for example asymptotics and congruences for $p(n)$. In a series of papers the author and Ono \cite{BO1,BO2} connected  the rank, a partition statistic introduced by Dyson, to weak Maass forms, 
 a new class of functions which are related to modular forms  and which were 
first considered in \cite{BF}. 
Here we do a further step towards understanding how weak Maass forms arise from interesting partition 
statistics by placing  certain  $2$-marked Durfee symbols introduced by Andrews \cite{An1}  into 
  the framework of weak Maass forms. 
  To do this we construct a new class of functions which we call quasiweak Maass forms 
  because they have quasimodular forms  as  components.  As an application we prove two conjectures of Andrews. 
  It seems that this new class of functions will play an important role in 
better understanding weak Maass forms of higher weight themselves, and  also their derivatives. 
  As a side product we introduce a new method which enables us   to prove transformation laws for generating functions over incomplete lattices.
\end{abstract}
\section{Introduction and Statement of Results}\label{intro} 
A partition of a  nonnegative integer $n$ is a non-increasing 
sequence of positive integers whose sum is $n$.
Let $p(n)$ denote the number of partitions of $n$.  
We have the generating function
\begin{eqnarray} \label{partgen}
P(q)=P(z):= \sum_{n=0}^{\infty} p(n)\, q^{24n-1} 
= \frac{1}{\eta(24z)},
\end{eqnarray}
where   $q:=e^{2 \pi i z}$ and 
$\eta(z):=
q^{ \frac{1}{24}}\, \prod_{n=1}^{\infty}(1-q^n)$
 is Dedekind's $\eta$-function. Of the many consequences of the modularity of $P(q)$, two of the most striking  are Rademacher's exact  formula for $p(n)$,  and Ramanujan type congruences. 
 In order to state the exact formula,  let $I_{s}(x)$ be the usual $I$-Bessel
function of order  $s$. Furthermore,
if $k$  and $n$ are positive integers, then define the Kloosterman sum 
\begin{equation*}
A_k(n):=
\sum_{ h \pmod k^*} \omega_{h,k}\,e^{ -\frac{2 \pi i hn}{k}},
\end{equation*}
where  the sum only runs over those $h$ modulo $k$ that are coprime to $k$,  and where
\begin{eqnarray*}
\omega_{h,k} :=
\exp\left(\pi i s(h,k) \right),
\end{eqnarray*}
with 
\begin{eqnarray*}
s(h,k):= \sum_{\mu \pmod k}  \left( \left( \frac{\mu}{k}\right) \right)   \left( \left( \frac{h \mu}{k}\right) \right).
\end{eqnarray*}
Here 
\begin{eqnarray*}
((x)):= \left \{ 
\begin{array}{ll}
x- \lfloor x \rfloor - \frac{1}{2} &\text{if } x \in \R \setminus \Z ,\\
0&\text{if } x \in \Z.
\end{array}
\right.
\end{eqnarray*}
If $n$ is a positive integer, then  Rademacher
showed that
\begin{equation}\label{Radformula}
p(n)= \frac{2 \pi}{ (24n-1)^{ \frac34}}
\sum_{k =1}^{\infty}
\frac{A_k(n)}{k}\cdot  I_{\frac{3}{2}}\left( \frac{\pi
\sqrt{24n-1}}{6k}\right).
\end{equation}
 In particular, (\ref{Radformula}) implies that 
 \begin{eqnarray*}
p(n)\sim \frac{1}{4n\sqrt{3}}\cdot e^{\pi\sqrt{2n/3}} \qquad  \qquad \text{as } n \to \infty.
\end{eqnarray*} 

The partition function    also has     nice congruence properties. For example,
   by  Ramanujan we have  that 
 \begin{eqnarray*}
p(5n+4) \equiv 0 \pmod 5,\\
p(7n+5) \equiv 0 \pmod 7, \\ 
p(11n+6) \equiv 0 \pmod {11},
\end{eqnarray*} 
for every $n \geq 0$ and   Ono \cite{On1}  showed that for any prime $\ell \geq 5$ there exist infinitely many non-nested arithmetic progressions  of the form $An+B$  such that
\begin{eqnarray*}
p(An+B) \equiv 0 \pmod \ell.
\end{eqnarray*}
In this paper we also consider asymptotics and congruences for certain partition  statistics.

To explain the Ramanujan congruences with modulus $5$ and $7$, Dyson \cite{Dy} introduced the  \textit{rank} of a partition, which  is defined to be its largest part minus the number of its parts. 
Dyson conjectured that the partitions of $5n+4$ (resp. $7n+5$) form $5$ (resp. $7$) groups of equal size when sorted by their ranks modulo $5$ (resp. $7$), and  this conjecture was proven by Atkin and Swinnerton-Dyer \cite{AS}.  
If $N(m,n)$ denotes the number of partitions of $n$  with rank  equal to $m$, then 
we have the generating  function 
\begin{equation} \label{RankGen}
R(w;q):= 1 + \sum_{m \in \Z} 
\sum_{n=1}^{\infty} 
N(m,n) w^m\, q^n
= 1+ \sum_{n=1}^{\infty} \frac{q^{n^2}}{(wq;q)_n
(w^{-1}q;q)_n}
= \frac{(1-w)}{(q;q)_{\infty}} \sum_{n \in \Z}
\frac{(-1)^n q^{\frac{n}{2}(3n+1) }}{1-wq^n},
\end{equation}
where 
$(a;q)_n:= \prod_{j=0}^{n-1} \left(1-aq^j \right)$ and  
$(a;q)_{\infty}:=\lim_{n \to \infty} (a;q)_n$.  
In particular,
\begin{eqnarray*}
R(1;q)&=&P(q),\\
R(-1;q)&=&f(q):=1+\sum_{n=1}^{\infty}\frac{q^{n^2}}{ (q;q)_n^2}.
\end{eqnarray*}
The function $f(q)$ is one of the  so-called \textit{mock theta functions}.  
Ramanujan listed 17 such functions in his last letter to Hardy and gave two more   in his ``Lost Notebook" \cite{Ra}, while  Watson \cite{Wa} defined three further functions.
Surprisingly, much remained unknown about these  series until recently. 
For example
 there was even debate concerning the rigorous definition of
such a function. 
Despite these  issues,
Ramanujan's mock theta functions have been shown to  possess many striking
properties, and they have been the subject of an astonishing  
number of important works 
(for example, see \cite{An1,
An3, AD, Ch, Co,
Dr,  Hi,  Zw}).
Much of this activity was foreshadowed by Dyson: 
 \\ 
 \\ 
\noindent
 ``\textit{The mock theta-functions give us tantalizing hints of a grand
 synthesis still to be discovered. Somehow it should be possible
 to build them into a coherent group-theoretical structure,
 analogous to the structure of modular forms which Hecke built
 around the old theta-functions of Jacobi. This remains a challenge for the
 future."}

\smallskip
\hskip3in Freeman Dyson, 1987

\hskip3in Ramanujan Centenary Conference
\\ \\ 
Recently much light has been shed on the nature of Ramanujan's mock theta functions.
Building on important work of Zwegers \cite{Zw}  the author and Ono (\cite{BO1,BO2})  solved  Dyson's ``challenge for  the future'' by placing the rank generating functions, and thus the mock theta function $f(q)$, in the context of weak Maass forms (see Section  \ref{CongSection}).
Loosely speaking it turns out that  for each root of unity   $w\not=1$, the function $R(w;q)$ is the  holomorphic part of a weight $\frac12$ weak Maass
form.   Since the  case $w=1$ yields a weight $-\frac12$ modular form, one can also 
view    the infinite family $R(w;q)$   as   a \textit{deformation} of $P(q)$.  

Viewing the rank generating functions in the framework of weak Maass forms has  
 found many applications, including  an exact formula
for the coefficients of $f(q)$ \cite{BO1}, asymptotics for    $N(m,n)$ \cite{Br}, and 
identities for rank differences \cite{BOR}. Moreover we obtain 
 congruences for $N(s,t;n)$ \cite{BO2}, the number of partitions of $n$ with rank congruent  to 
 $s$ modulo  $t$, which  give a combinatorial  decomposition of  congruences for $p(n)$.

Naturally it is of wide interest to find other explicit
examples of Maass forms. Here we construct a new infinite family
of such forms, arising in a more complicated way from   an interesting partition statistic introduced by  Andrews  \cite{An2}.  
Define the    \textit{symmetrized second moment function}   
\begin{eqnarray*}
\eta_2(n):= \sum_{m = - \infty}^{\infty} 
\left(
\begin{matrix}
m \\
2
\end{matrix}
\right)
N(m,n).
\end{eqnarray*}  
Andrews  showed that this   function  enumerates the 2-marked Durfee symbols (see  Section  \ref{CombiSection}). 
 We have the generating function 
\begin{eqnarray} \label{MomGen}
R_ 2(q):=
\sum_{n=0}^{\infty} \eta_{2}(n)\,q^n = \frac{1}{(q;q)_{\infty}} 
\sum_{n \in \Z \setminus \{ 0\} }
(-1)^{n+1} 
\frac{q^{\frac{3n}{2} (n+1)  }}{(1-q^n)^{2}}.
\end{eqnarray} 
One of the main goals of this paper is to place this function into the framework of weak Maass forms and 
introduce a new class of functions, which we call quasiweak Maass forms, which  seem to play an important role in  better understanding weak Maass forms themselves, and  also their derivatives. 
It will turn out that  the  function (\ref{MomGen}) is part of a weight $\frac32$ weak Maass forms and 
 occurs  from deforming the weight $\frac12$ rank generating function by applying a certain differential operator which is responsible for the higher weight.  
 The proof of the modularity of the usual rank generating function (\ref{RankGen}) heavily relies on the fact that  the sum runs through the full lattice $\Z$ which allows Poisson summation. 
In contrast, the sum of the generating function  in (\ref{MomGen}) runs only  through the  incomplete lattice  
$\Z \setminus \{ 0\}$.  To our knowledge so far known examples which use Poisson summation to prove
modularity come from sums over full lattices.
The first guess to deal with missing terms in Poisson summation seems to be to add extra summands
during the computations, and remove them at the end. 
However, this does not work in the present case due to the double pole in the $n=0$ term. 

We overcome those problems by introducing a new method which enables us to handle incomplete sums like (\ref{MomGen}). 
For this purpose, we embed (\ref{MomGen}) into a bigger family   which includes $R_2(q)$ as a limiting case of a certain differential operator. 
It seems very likely that this new method will be helpful in other  settings as well, and may lead to a better understanding of functions like $R(w;q)/(1-w)$ (a modified version of Dyson's rank generating function that removes the ``artificial'' zero).  In our situation a very careful analysis of introduced poles is required since they give rise to additional terms in the transformation law, which we managed to identify as quasimodular components.   
Recall that a meromorphic function $f:\H \to \C$  is called a \textit{quasimodular form} if 
it   can be  written as linear combination of derivatives of modular forms.
Note that some authors do not include meromorphic functions in this definition and only allow forms that lie in the algebra generated by $E_2, E_4,$ and $E_6$.

To state our result   define 
\begin{eqnarray*}
\mathcal{R}(z) := R_2 (24 z) e^{ - 2\pi i z},
\end{eqnarray*}
\begin{eqnarray}\label{Mordelint}
\mathcal{N}(z):= 
\frac{i}{4\sqrt{2}\pi} 
\int_{- \bar z}^{i \infty} 
\frac{\eta(24\tau )}{(- i (\tau+z))^{\frac32}} \, d\tau,
\end{eqnarray}
and 
\begin{eqnarray}
\mathcal{M} (z) :=
 \mathcal{R}(z)- 
 \mathcal{N}(z)  - 
 \frac{1}{24\eta(24 z)} 
 +
 \frac{E_2(24 z)}{8\eta(24 z)},
\end{eqnarray}
where  as usual 
$$E_2(z):=  
1- 24 \sum_{n=1}^{\infty} \sigma_1(n)\,q^n
$$ 
with $\sigma_1(n) := \sum_{d|n}d$. 
The function  $E_2(z)$  is ``nearly modular'', and in particular      satisfies
$E_2(z+1)=E_2(z)$ and 
\begin{eqnarray} \label{E2transf}
E_2\left( -\frac{1}{z} \right) 
= z^2E_2(z) +\frac{6  z}{\pi i} .
\end{eqnarray} 
\begin{theorem} \label{maintheorem2}
The function $\mathcal{M}(z)$ is a harmonic weak Maass form of weight $\frac{3}{2}$ 
on  $\Gamma_0(576)$ with Nebentypus character $\chi_{12}(\cdot):=\leg{12}{\cdot}$.
\end{theorem}
\textit{Two remarks.}

\noindent 1) 
Since $\mathcal{R}(z)$ is, up to quasimodular forms, the holomorphic part of a weak Maass form, we refer  to $R_2(q)$ as a \textit{quasimock theta function}.

\noindent 2)
The function $R_2(q)$ is also  crucial  for understanding $k$-marked Durfee symbols for $k >2$ (see \cite{An2} for the definition). In upcoming work Mahlburg  and the author show  \cite{BM} that modularity properties of the generating functions $R_k(q)$ can  be concluded from modularity properties of $R(q)$ and $R_2(q)$, and that these functions give rise to a new class of functions which we call \textit{quasiweak Maass forms}.  It seems that this new class of functions will play an important role in 
better understanding weak Maass forms of higher weight themselves, and  also their derivatives.

The fact that  $\mathcal{M}(z)$ is a weak Maass form has some nice applications.
Here we address two of them: congruences  and asymptotics.
These were formulated by Andrews as open problems
(problems 11 and 13 page 39 of \cite{An2}).
\begin{theorem} \label{maintheorem1}
We have  
\begin{multline*}
\eta_2(n) = \sum_{k=1}^{\left[ \sqrt{n}\right] } A_k(n)
\left(     - 
\frac{3}{2
  (24n-1)^{\frac{1}{4}}}
 I_{\frac{1}{2}} \left( \frac{\pi}{6k} \sqrt{ 24n-1}   \right) 
+ \frac{\pi(24n-1)^{\frac{1}{4}} }{12k} \right. \\ 
\left. \, I_{-\frac{1}{2}} \left(\frac{\pi}{6k} \sqrt{24n-1} \right)
+  \frac{\pi}{12k (24n-1)^{\frac{3}{4}}} \, I_{\frac{3}{2}} \left(\frac{\pi}{6k} \sqrt{24n-1} \right)
\right) 
+ O \left( n^{1+ \epsilon }\right).
\end{multline*}
\end{theorem}
\noindent
In particular the $k=1$ term gives the main term in the asymptotics expansion.
\begin{corollary} \label{corasymp1}
As $n \to \infty$
\begin{eqnarray*}
\eta_2(n) \sim \frac{1}{4 \sqrt{3}} \, e^{ \frac{\pi}{6} \sqrt{24n-1}}.
\end{eqnarray*}
\end{corollary}  
\begin{remark}
In particular  Corollary \ref{corasymp1} implies that as $n \to \infty$
$$
\eta_2(n) \sim n p(n).
$$
\end{remark}
We next turn to congruences.
As in the case of partitions one can associate certain ranks to $2$-marked Durfee symbols, which we will describe in Section \ref{CombiSection}. 
Let $NF_2(r,t;n)$ denote the number of $2$-marked Durfee symbols with (full) rank congruent to $r$ modulo $t$.   
\begin{theorem}\label{congtheorem}
Let $t>3$ be an   integer, $j \in \N$,  and  $\mathcal{Q} \nmid 6t$  a prime. Then there exist infinitely many arithmetic progressions $An+B$ such that for every $ 0 \leq r<t$, we have 
\begin{displaymath}
NF_2(r,t;An+B) \equiv 0 \pmod{\mathcal{Q}^j}.
\end{displaymath}
\end{theorem}   
\begin{remark}
With similar methods as in   \cite{B2},  one could also obtain congruences for $t=\mathcal{Q}^{\ell}$ with $\ell \in \N$.
\end{remark}
Theorem \ref{congtheorem}  can be viewed as yielding congruences for $\eta_2(n)$ as well as a combinatorial decomposition of these congruences.
\begin{corollary} \label{congcorollary}
Let $j$ be a positive  integer  and $\mathcal{Q}>3$ a prime. 
Then there exist  infinitely many arithmetic progressions $An+B$ such that   \begin{displaymath}
\eta_2(An+B) \equiv 0 \pmod{\mathcal{Q}^j}.
\end{displaymath}
\end{corollary}
From Theorem \ref{maintheorem1} and Corollary \ref{congcorollary} one can obtain  congruences and asymptotics  for an interesting new  statistic introduced by Andrews \cite{An4}.
Let $spt(n)$ be the total number of appearances of smallest parts in each integer partition of $n$.  
Andrews showed that 
$$
spt(n) = n p(n) - \eta_2(n).
$$
Since congruences and asymptotics are known for $p(n)$, we should also be able to obtain similar 
results for $spt(n)$. However, we do not further address this topic here. 

We next consider odd moments
\begin{eqnarray*}
\eta_2^o(n):= 
\sum_{m \in \Z} 
\left( \begin{matrix} m+1 \\ 
2 \end{matrix}\right)
N^o(m,n),
\end{eqnarray*}
where $N^o(m,n)$ is the number of partitions related to an odd Durfee symbol with odd rank $m$ (for the definition of odd ranks and odd Durfee symbols see Section  \ref{CombiSection}).
We have    
\begin{eqnarray*}
R_{2}^o(q):=
\sum_{n=1}^{\infty}
\eta_{2}^o(n)\, q^n 
= \frac{1}{(q^2;q^2)_{\infty}}
 \sum_{n \in \Z}
 (-1)^n\frac{q^{3n^2+5n + 2}}{\left( 1-q^{2n+1}\right)^{3}}.
 \end{eqnarray*}  
 Defining  the  sum of Kloosterman type 
\begin{eqnarray*}
A^o_k(n):= e^{\frac{\pi i k}{2} }\sum_{h \pmod k^*} 
\omega_{2h,k}\, e^{\frac{3 \pi i hk}{2} }\, e^{\frac{\pi i h}{2k} (1-4n)}
\end{eqnarray*}
we find the following asymptotic expansion for $\eta_2^o(n)$.
\begin{theorem} \label{maintheorem3}
We have 
\begin{multline*}
\eta_2^o(n)
= -\frac{i}{\sqrt{2}} \sum_{k \text{ odd} }^{ \left[ \sqrt{n}\right]}
A_k^o(n)
\left(
-\frac{3\pi^2(3n-1)^{\frac{1}{4}}}{4}  I_{-\frac{1}{2} } \left( \frac{\pi}{3k} \sqrt{3n-1}\right) \right. 
\\
\left.
+ \frac{\pi (3n-1)^{\frac{3}{4}}}{16k}  I_{-\frac{3}{2} } \left( \frac{\pi}{3k} \sqrt{3n-1}\right)
\right) + O(n^{2+\epsilon}).
\end{multline*}  
 \end{theorem} 
 \begin{corollary}
 As $n \to \infty$,
 \begin{displaymath}
 \eta_2^o(n) \sim 
 \frac{3}{8} \sqrt{n} \, e^{ \frac{\pi}{3} \sqrt{3n-1}}.
 \end{displaymath}
 \end{corollary}
 The proof of Theorem \ref{maintheorem3} heavily relies on  a transformation law for $R_2^o(q)$ which is more complicated than the transformation law for $R_2(q)$ arising from higher  order poles of the generating function. For this reason no   usual weak Maass forms arise. 
 It would be very interesting to analyze the class of functions $R_2^o(q)$ belong to.  
 
 The paper is organized as follows. 
 In  Section \ref{CombiSection}  we recall the connection between higher moments and marked Durfee symbols.
 In Section \ref{TRANS} we prove a transformation  law for $R_2(q)$. To overcome the above mentioned problems,  we introduce a new family of functions  
 $\mathcal{R}_2(w;q)$
 which have $R_2(q)$ as a limiting case of a certain differential operator $L$.  
We  first prove  a transformation law for $\mathcal{R}_2(w;q)$, and then apply $L$. 
Here a careful analysis is required since we introduce summands which have poles when  $w=1$, and it is because of these terms that  quasimodular forms show up.
 In Section  \ref{MaassSection} we realize  the  error integral that arises as an integral of theta functions and prove Theorem \ref{maintheorem1}.
 As an application we then prove asymptotics in Section \ref{AsympSection} and congruences in Section 
 \ref{CongSection}.  In Sections \ref{TransoddSection} and \ref{AsympoddSection}  we consider odd rank generating functions. 
 \section*{Acknowledgements}
\noindent
The author thanks Frank Garvan, Karl Mahlburg, David Penniston, Rob Rhoades,  and 
the referee for helpful comments.  
 \section{Some combinatorial results on $2$-marked Durfee symbols} \label{CombiSection}
 Here we recall some results on $2$-marked Durfee symbols and their connection to rank moments \cite{An2}.
 Recall that the largest square of nodes in the Ferrers graph of a partition is called the \textit{Durfee square}.  
The \textit{Durfee symbol} consists of 2 rows and a subscript, where the top row consists of the columns to the right of  the Durfee square,   the bottom row consists of the rows  below the Durfee square, and the subscript denotes the side length of the Durfee square.
The number being partitioned is equal to the sum of the rows of the symbol   plus the number of nodes in the Durfee  square. 
For example the Durfee symbol
\begin{eqnarray*}
\left(
\begin{matrix}
2& \\
3&1
\end{matrix}
\right)_4
\end{eqnarray*}
represents a partition of $2+3+1+4^2=22$.
To define \textit{2-marked Durfee symbols} one requires two copies of  the positive integers denoted by $\{ 1_1,2_1,\cdots \}$ and $\{ 1_2,2_2,\cdots \}.$
We   form Durfee symbols as before and use the two copies  of the positive  integers  for parts in both rows.
We additionally demand that the following conditions are met:
\begin{enumerate}
\item
The sequence of parts and the sequence of subscripts in each row are non-increasing.
\item The subscript $1$ occurs at least once in the top row.
\item If $M_1$ is the largest part with  subscript $1$  in the top row, then all parts in the bottom row with subscript 1 lie in $[1,M_1]$ and with subscript $2$ lie in $[M_1,S]$, where $S$ is the side of the Durfee symbol. 
\end{enumerate}
We  denote by $\mathcal{D}_2(n)$ the number of $2$-marked Durfee  symbols arising  from partitions of $n$.
In \cite{An2} it is shown that 
\begin{eqnarray*}
\mathcal{D}_{2} (n) = \eta_{2}(n).
\end{eqnarray*} 
As in the case of partitions one can associate ranks to  $2$-marked Durfee symbols \cite{An2}. For a $2$-marked Durfee symbol $\delta$ we define the \textit{full rank} $FR(\delta)$ by 
\begin{eqnarray*}
FR(\delta) := \rho_1(\delta) + 2 \rho_2(\delta),
\end{eqnarray*}
where
\begin{eqnarray*}
\rho_i(\delta) := \left\{ 
\begin{array}{ll}
\tau_i(\delta) - \beta_i(\delta)-1 & \text{ for }i=1,\\
\tau_i(\delta)- \beta_i(\delta) &\text{ for } i=2.
\end{array}
\right. 
\end{eqnarray*}
Here $\tau_i(\delta)$ (resp. $\beta_i(\delta)$) denotes the number of entries in the top (resp. bottom) row of $\delta$ with subscript $i$.
We let $NF_2(m,n)$ denote the number of 2-marked Durfee symbols related to $n$ 
with full rank $m$ and $NF_2(r,t;n)$ the number of $2$-marked Durfee symbols related to $n$ with full rank congruent to $r$  modulo  $t$. 
Let 
\begin{eqnarray} \label{2generating}
R_2(w;q):=
\sum_{n=1}^{\infty} \sum_{m \in \Z} NF_2(m,n)\, w^m\, q^n.
\end{eqnarray} 
In particular  we have 
\begin{eqnarray*}
R_2(1;q)=R_2(q).
\end{eqnarray*}
Moreover if $w^3 \not =1$, then  
$$
R_2(w;q)= \frac{w^2}{(1-w)(w^3-1)}
\left(R(w;q)- R(w^2;q) \right).
$$   
Therefore in this case   properties of  $R_2$, like congruences or asymptotics,  can be concluded from properties of the usual rank generating function $R$ using work of the author and Ono.

We next consider  odd Durfee symbols.  
 A partition into pairs of consecutive integers (excepting one instance of the largest part) has associated to it the \textit{odd Durfee symbol}
$
\left(  \begin{smallmatrix} a_1&a_2&\cdots&a_i \\
b_1&b_2&\cdots&b_j
\end{smallmatrix}\right)_n 
$
with  the $a's$ and $b's$ being  odd numbers $\leq 2n+1$,  
 and the number being partitioned is 
$2n^2+2n+1+ \sum_{m=1}^{i} a_m+ \sum_{k=1}^{j} b_k$.
The \textit{odd rank} of an odd Durfee symbol is defined as  the number of entries in the top minus the number of entries in the bottom row.
We denote by $N^o(m,n)$ the number of partitions related to an odd Durfee symbol with odd rank $m$.
We have the generating function
\begin{eqnarray*}
R^o(w;q)
:= \sum_{n=1}^{\infty}
\sum_{m \in \Z} N^o(m,n) w^m q^n 
= \frac{1}{\left(q^2;q^2  \right)_{\infty}} 
\sum_{n \in \Z} (-1)^n 
\frac{q^{ 3n^2+3n+1}}{1-wq^{2n+1}}.
\end{eqnarray*}
The definition of an  \textit{odd 2-marked Durfee symbol}   is almost the same as for a 2-marked Durfee symbol. 
The only changes are:
\begin{enumerate}
\item
All entries in the symbol are odd numbers.
\item 
The subscribing rules are the same but instead of adding $n^2$ to the sum
 when the Durfee symbol is of side $n$ we now add $2n^2+2n+1$. 
\end{enumerate}
 If  $\mathcal{D}_2^o(n)$ denotes the number of odd $2$-marked Durfee symbols of $n$, then \cite{An2}
 \begin{eqnarray*}
 \mathcal{D}_{2}^o(n) = \eta_{2}^o(n).
 \end{eqnarray*}   

\section{A transformation law  for $R_2(q)$} \label{TRANS}  
Due to  problems   mentioned in  the introduction  we cannot show a transformation law for  $R_2(q)$ directly but we introduce a function of an additional variable $w$ that is related to $R_2(q)$.
For this purpose, define for 
$w \in \C$ with Re$(w)$ sufficiently small
\begin{eqnarray} \label{wgen}
\mathcal{R}_{2}(w;q)=
\mathcal{R}_{2}(w;\tau):=
- \frac{1}{(q;q)_{\infty}}
\sum_{n \in \Z \setminus \{0 \}}  
\frac{(-1)^{n} q^{ \frac{n(3n+1)}{2} }}{1-e^{2 \pi i w}\, q^n},
\end{eqnarray}
where $q:=e^{2 \pi i \tau }$. Defining the operator $L$ for  a function $g$ that is differentiable in a neighborhood of $0$,  
$$
L(g(w)) :=   \frac{1}{2 \pi i } 
\left[\frac{\partial}{\partial w} g \right]_{w=0},
$$
we can connect the function $R_2(q;w)$ to $R_2(q)$ 
$$
R_2(q) =  L(\mathcal{R}_{2}(w;q)).
$$
We prove a transformation law for $\mathcal{R}_2(w;q)$ from which  we  conclude a  transformation law for $R_{2}(q)$ by applying $L$. 
For this   let
$$ 
H_w(x):= \frac{e^x}{1-e^{2 \pi iw}\, e^{2x}}.
$$
Then 
$$
H_w(x+  \pi  i) = - H_w(x).
$$
Moreover   for an integer $\nu$ and $z \in \C$ with Re$(z)>0$, define 
\begin{eqnarray*}
I_{k,\nu}^{\pm} (z;w) :=
 \int_{\R} 
 e^{ -\frac{3 \pi z x^2}{k} } 
  \, H_w \left( \pm \left( \frac{\pi i \nu}{k}  -\frac{\pi i}{6k}  -\frac{\pi z x}{k}\right)  \right)\, dx.
  \end{eqnarray*}
\begin{theorem} \label{transtheorem1}
For coprime integers $h$ and $k$ with $k>0$, let $q:=e^{\frac{2 \pi i }{k} (h+iz)}$ and     $q_1:=e^{\frac{2 \pi i }{k} \left(h'+\frac{i}{z}\right)}$, where for $k$ odd we let $h'$ be an even solution of  $hh' \equiv - 1\pmod k$, and for  even $k$ we let  $h'$ be defined by $hh'\equiv -1 \pmod{2k}$.
 Then 
 \begin{multline*}
 \mathcal{R}_{2}(w;q) 
= \frac{\omega_{h,k} \, z^{\frac{1}{2}} \, e^{\frac{\pi}{12k}\left(z^{-1} -z\right)}}{(q_1;q_1)_{\infty}}
 \left( 
 \frac{1}{1-e^{2 \pi i w}} - \frac{i\, e^{\frac{3 \pi k w^2}{z}- \pi i w -\frac{\pi w}{z}} }{z \left(1-e^{-\frac{2\pi w}{z}}  \right)}  
 \right) \\   +
 i z^{-\frac{1}{2}} \omega_{h,k} e^{ \frac{\pi}{12 k}\left( z^{-1} - z\right) } 
  e^{\frac{3 \pi kw^2}{z} - \pi i w - \frac{\pi w}{z}}\, \mathcal{R}_2\left(\frac{iw}{z} ;q_1\right)  \\
 -  \frac{z^{\frac{1}{2}}}{k} \, \omega_{h,k} \, e^{-\frac{\pi z}{12k} } 
 \sum_{\substack{\nu \pmod k \\ \pm  }}  
(-1)^{\nu} \, 
 e^{  \frac{\pi i h'}{k}(-3 \nu^2 + \nu) }  
 I_{k,\nu}^{\pm}(z;w).
 \end{multline*}
\end{theorem}
\begin{proof} 
We first observe that the difference of both sides of the transformation law  is a holomorphic function of $w$ in some neighborhood of $0$.
Indeed this follows immediately if $w\not=0$ is sufficiently small.
Moreover it  follows from (\ref{hol}) that   the point $w=0$ is a removable singularity. 
Therefore it is  by  the identity theorem of holomorphic  functions  enough to prove the theorem   for $w \in \C$  with Re$(w)\not =0$ sufficiently small and $z>0$ real.  For those $w$, the function  
\begin{eqnarray*}
\widetilde{R}_2(w;q):= 
\sum_{n \in \Z}  
\frac{(-1)^{n-1} q^{ \frac{n(3n+1)}{2}  }}{1-e^{2 \pi i w}\, q^n}
\end{eqnarray*}
is a holomorphic  function of $z$.
Writing $n= \nu + km$ with $m \in \Z$ and $\nu$ running modulo $k$, we obtain
\begin{equation} \label{poisson2}
\widetilde{R}_{2}(w;q) 
= - \sum_{\nu \pmod k} 
(-1)^{\nu} \, e^{ \frac{3 \pi i h \nu^2}{k} }  
\sum_{ m \in \Z} (-1)^m 
e^{- \frac{3\pi z}{k} ( \nu + km)^2   }\,
H_w \left(  \frac{\pi i \nu h}{k} - \frac{\pi z}{k} (\nu+ km)\right)\, 
.
\end{equation}
Poisson summation gives that the inner sum equals
\begin{eqnarray*} &&
\sum_{ n \in \Z}  
\int_{\R} 
e^{ - \frac{3 \pi z}{k} (\nu + kx)^2 } 
H_w \left(  \frac{\pi i \nu h}{k} - \frac{\pi z}{k} (\nu+ kx)\right)\, e^{ 
 \pi i (2n+1)x} \, dx 
\\ 
&=&  
\frac{1}{k}
\sum_{ n \in \Z} 
\int_{\R}
e^{ - \frac{3 \pi zx^2}{k}  } 
H_w \left(  \frac{\pi i \nu h}{k} - \frac{\pi zx}{k}  \right)\, e^{ \frac{\pi i}{k} (2n+1)(x- \nu)} \, dx,
\end{eqnarray*}
where for the last equality we made the change of variables $x \mapsto \nu+ kx$.
Inserting this  into (\ref{poisson2}), and changing  in 
 the sum over $n \leq -1$,  $n \mapsto -n-1$, $x \mapsto -x$, and $\nu \mapsto - \nu$, gives  \begin{equation*}
\widetilde{R}_{2}(w;q)=-
\frac{1}{k} 
\sum_{\substack{\nu \pmod k    \\   \pm  }  }
(-1)^{\nu} \, e^{ \frac{3 \pi i h \nu^2}{k}   }  \\
\sum_{n \geq 0}
\int_{\R} 
e^{ - \frac{3 \pi zx^2}{k}   } 
H_w \left( \pm \left(  \frac{\pi i \nu h}{k} - \frac{\pi z x}{k}\right)\right)\, e^{ \frac{\pi i}{k} (2n+1)(x- \nu)} \, dx.
\end{equation*}
Next  write  
\begin{eqnarray*}
H_w \left( \pm \left(  \frac{\pi i \nu h}{k} - \frac{\pi zx}{k}  \right)\right) 
= \frac{\pm \, e^{- \pi i w  }}{ 2 \sinh \left(  \frac{\pi z x}{k} - \frac{\pi i h \nu}{k}  \mp  \pi i w\right)}.
\end{eqnarray*} 
Using that the function 
\begin{eqnarray} \label{Sw}
S_w^{\pm}(x) 
:= \frac{\sinh ( x \pm  k \pi  i w)}{\sinh \left(\frac{x}{k}  \pm  \pi i w  \right)}
\end{eqnarray}
is an entire function,  we see   that the only possible poles of the integrand lie in points
\begin{equation} \label{xm}
x_m^{\pm} 
:= \frac{i}{z} \left(m \pm k w\right).
\end{equation}
Moreover for a fixed $m$ this leads at most    to a pole if 
$h \nu \equiv m \pmod  k$,
 and we may choose  
\begin{equation} \label{nm}
\nu_m:=  - h' m.
\end{equation}
Thus poles may at most occur for $x_m$ and $\nu_m$ as in (\ref{xm}) and (\ref{nm}).
Using the Residue Theorem we  shift the path of integration through 
$$ 
\omega_n := \frac{(2n+1)i}{6z}.
$$
For this  we have to take those $x_m^{\pm}$ into account for which $n \geq 3m\geq \frac12 ( 1 \mp 1)$. We denote the associated residues of each summand by $\lambda_{n,m}^{\pm}$.
Then
$$
\widetilde{R}_{2}(w;q) = \sum_1 + \sum_2,
$$
where 
\begin{eqnarray*}
\sum_1 &:=& - \frac{2 \pi i}{k} 
\left( 
\sum_{m \geq 0} \sum_{n \geq 3m} \lambda_{n,m}^+ 
+  \sum_{m > 0} \sum_{n \geq 3m} \lambda_{n,m}^-
\right),
\\
\sum_2 &:= &
- 
\frac{1}{k} 
\sum_{\substack{\nu \pmod k    \\  \pm  }  }
(-1)^{\nu} \, e^{ \frac{3 \pi i h \nu^2}{k}} 
\sum_{n \geq 0} \int_{\R+ \omega_n} 
e^{ - \frac{3 \pi zx^2}{k}   } 
H_w \left( \pm  \left(  \frac{\pi i \nu h}{k} - \frac{\pi zx}{k} \right)\right)\, e^{ \frac{\pi i}{k} (2n+1)(x- \nu)} \, dx.
\end{eqnarray*}
We first consider $\sum_1$.
We have
$$
\lambda_{n,m}^{\pm} 
= \pm k  (-1)^{\nu_m}\, 
\frac{e^{- \frac{3\pi z \left.x_m^{\pm}\right.^2}{k}  + \frac{\pi i}{k}(2n+1)(x_m^{\pm} 
- \nu_m) + 
\frac{3\pi i h \nu_m^2}{k}  - \pi i w 
 }}{  2 \cosh \left(\frac{\pi z x_m^{\pm}}{k}  - \frac{\pi i h \nu_m}{k }  \mp \pi i w\right) \pi z},
$$
which yields  
$$
\lambda_{n+1,m}^{\pm} = e^{\frac{2 \pi i}{k}  (x_m^{\pm}- \nu_m)  }  \lambda_{n,m}^{\pm} \,.
$$
Therefore 
$$
\sum_{n =3m}^{\infty} \lambda_{n,m}^{\pm} 
= \frac{\lambda_{3m,m}^{\pm}}{1- \exp \left(\frac{2 \pi i}{k} (x_m^{\pm} -\nu_m)  \right)}.
$$
Using this one can prove that 
\begin{eqnarray*}
\sum_1 
=  \frac{i}{z}\, e^{\frac{3 \pi kw^2}{z}  - \pi iw - \frac{\pi w}{z}} 
\widetilde{R}_2\left(q_1;\frac{iw}{z} \right)  .
 \end{eqnarray*}
 We next turn to $\sum_2$. 
 With the same argument as before, we   change the sum  over $n \in \N_0$ back into a sum over all $n \in \Z$.
 Substituting $x \mapsto x + \omega_n$ and writing $n=3p+ \delta$ with $p \in \Z$ and $\delta \in \{ 0,\pm1\}$ gives 
 \begin{multline*}
 \sum_2 = 
 - \frac{1}{k} \sum_{ \nu,p,\delta  } (-1)^{\nu} \, e^{ \frac{3 \pi i h \nu^2}{k} - \frac{\pi(6p+2 \delta+1)^2}{12kz} - \frac{\pi i}{k}(6p+ 2 \delta +1) \nu } \\
 \int_{\R} 
 e^{- \frac{3 \pi z x^2}{k}}
 H_w \left(  
 \frac{\pi i h \nu}{k} - \frac{\pi z x}{k} - \frac{\pi i (6p+2 \delta+1)}{6k}
 \right)\, dx.
 \end{multline*} 
 Changing $ \nu$ into $- h'(\nu+p)$, a lengthy calculation gives
 \begin{multline*}
 \sum_2 = - \frac{1}{k} 
 \sum_{ \nu, \delta,p  } 
 (-1)^{\nu+p} 
 e^{ \frac{\pi i h'}{k}(-3 \nu^2+ (2 \delta+1) \nu) - \frac{\pi }{12kz}(2 \delta+1)^2} \\
 q_1^{\frac{p}{2} (3 p + 2 \delta+1) }
 \int_{\R} 
 e^{ -\frac{3 \pi z x^2}{k} } H_w \left(\frac{\pi i \nu}{k} - \frac{\pi i}{6k}(2 \delta+1) -\frac{\pi zx}{k} \right)\, dx.
 \end{multline*}
 Now the integral is independent of $p$ and the sum over $p$ equals 
 \begin{eqnarray} \label{eta}
 \sum_{p \in \Z} (-1)^p \, q_1^{ \frac{p}{2}(3p+2 \delta +1)}.
 \end{eqnarray}
 If $\delta=1$, then (\ref{eta}) vanishes since the $p$th and the $-(p+1)$th term cancel.
 If $\delta=0$ or  $\delta=-1$, then   (\ref{eta}) equals $(q_1;q_1)_{\infty}$.
 Changing for $\delta=-1$, $x \mapsto - x$ and $ \nu \mapsto - \nu$ gives
 \begin{eqnarray*}
 \sum_2  = - \frac{(q_1;q_1)_{\infty}}{k}  e^{- \frac{\pi}{12kz} }
 \sum_ {\substack{\nu \pmod k    \\  \pm  }   } (-1)^{\nu} 
  e^{\frac{\pi i h'}{k} (-3 \nu^2 + \nu  ) } 
 I_{k,\nu}^{\pm}(z;w).
  \end{eqnarray*}
  From this the theorem follows using
  \begin{eqnarray*}
  (q_1;q_1)_{\infty} = \omega_{h,k}\, z^{\frac{1}{2}} \, e^{\frac{\pi}{12k} (z^{-1}- z)} (q;q)_{\infty}.
  \end{eqnarray*}
\end{proof}
From Theorem \ref{transtheorem1} we  conclude the desired transformation law for $R_2(q)$.
\begin{corollary} \label{transcorollary1}
For coprime integers $h$ and $k$ with $k>0$, let $q:=e^{\frac{2 \pi i }{k} (h+iz)}$ and     $q_1:=e^{\frac{2 \pi i }{k} \left(h'+\frac{i}{z}\right)}$, where for $k$ odd we let $h'$ be an even solution of  $hh' \equiv - 1\pmod k$, and for  even $k$ we let  $h'$ be defined by $hh'\equiv -1 \pmod{2k}$.
 Then   
\begin{multline*}
R_2(q)
=  \omega_{h,k}\, z^{\frac{1}{2}}\, e^{\frac{\pi}{12k} (z^{-1}- z)}   
\left( 
\frac{1}{(q_1;q_1)_{\infty} } \left(
-\frac{3k}{4\pi z} + \frac{1}{24z^2} + \frac{1}{24} \right) 
+ \frac{1}{z^2} R_2(q_1)   \right) \\
- \frac{z^{\frac{1}{2}} }{2k}\, \omega_{h,k}\, e^{- \frac{\pi z}{12k}}
\sum_{\nu \pmod k   }   (-1)^{\nu}
e^{ \frac{\pi i h'}{k} (-3 \nu^2+ \nu)}
I_{k,\nu}(z),
\end{multline*}
where 
$$
I_{k,\nu}(z):=\int_{\R} e^{- \frac{3 \pi z x^2}{k}}\, 
\frac{\cosh \left(\frac{\pi i \nu}{k} - \frac{\pi i}{6k}  - \frac{\pi zx}{k} \right)}{\sinh^2 \left(\frac{\pi i \nu}{k} - \frac{\pi i}{6k}  - \frac{\pi zx}{k} \right)} \, dx.
$$
\end{corollary}
\begin{proof}
One can show that 
\begin{eqnarray}\label{hol}
L \left( \frac{1}{1-e^{2 \pi i w}  } - \frac{i\, e^{ 3\pi k w^2 - \pi i w  - \frac{\pi w}{z} }}{z\left( 1- e^{- \frac{2 \pi w}{z} } \right)}\right) &=& 
- \frac{3 k}{4 \pi z} + \frac{1}{24} + \frac{1}{24z^2},
\end{eqnarray}
\begin{eqnarray*}
L \left(    
   \frac{e^{ 3\pi k w^2 - \pi i w  - \frac{\pi w}{z} }}{1- q_1^m\,e^{- \frac{2 \pi w}{z} }}\right) &=& 
   - \frac{1}{2(1-q_1^m)} + \frac{i(1+q_1^m)}{2z(1-q_1^m)^2},
   \\ 
   L \left( H_w \left(\frac{\pi i \nu}{k} - \frac{\pi i}{6k}- \frac{\pi z x}{k} \right)  
   +  H_w \left(-\frac{\pi i \nu}{k} + \frac{\pi i}{6k} + \frac{\pi z x}{k} \right)
   \right)
   &=& \frac{\cosh \left( \frac{\pi i \nu}{k}  - \frac{\pi i}{6k} - \frac{\pi z x }{k} \right)}{2 \sinh^2 \left(\frac{\pi i  \nu}{k}   
   - \frac{\pi i}{6k} - \frac{\pi z x}{k}   \right)}.
 \end{eqnarray*}
 Moreover 
 \begin{eqnarray*}
 \sum_{m \in \Z \setminus \{ 0\} }
 \frac{(-1)^m q^{\frac{m}{2} (3m+1) }\left(1+q^m\right)}{\left( 1-q^m\right)^2} 
 &=& 2   \sum_{m \in \Z \setminus \{ 0\} }
 \frac{(-1)^m q^{\frac{m}{2} (3m+1) +m}  }{\left( 1-q^m\right)^2} 
 =  - 2 (q;q)_{\infty} R_2(q),\\
 \sum_{m \in \Z\setminus \{0 \}} 
 \frac{(-1)^m\, q^{\frac{m}{2}(3m+1) }}{1-q^m} &=& 0,
 \end{eqnarray*}
 since the $m$th and $-m$th term cancel. 
  Now the corollary follows from Theorem \ref{transtheorem1}.
\end{proof}
\section{Proof of Theorem \ref{maintheorem2}} \label{MaassSection}
In this section we prove that $\mathcal{M}(z)$ is a weak Maass form.  
If $k\in \frac{1}{2}\Z\setminus
\Z$, 
$z=x+iy$ with $x, y\in \R$, then the weight $k$ hyperbolic
Laplacian is given by
\begin{equation}\label{laplacian}
\Delta_k := -y^2\left( \frac{\partial^2}{\partial x^2} +
\frac{\partial^2}{\partial y^2}\right) + iky\left(
\frac{\partial}{\partial x}+i \frac{\partial}{\partial y}\right).
\end{equation}
If $v$ is odd, then define $\epsilon_v$ by
\begin{equation}
\epsilon_v:=\begin{cases} 1 \ \ \ \ &{\text {\rm if}}\ v\equiv
1\pmod 4,\\
i \ \ \ \ &{\text {\rm if}}\ v\equiv 3\pmod 4. \end{cases}
\end{equation}
Moreover we let $\chi$ be a Dirichlet character.
 A {\it (harmonic) weak Maass form of weight $k$ with Nebentypus $\chi$ on a subgroup
$\Gamma \subset \Gamma_0(4)$} is any smooth function $g:\H\to \C$
satisfying the following:
\begin{enumerate}
\item For all $A= \left(\begin{smallmatrix}a&b\\c&d
\end{smallmatrix} \right)\in \Gamma$ and all $z\in \H$, we
have 
\begin{displaymath}
g(Az)= \leg{c}{d}^{2k}\epsilon_d^{-2k} \chi(d)\,(cz+d)^{k}\ g(z).
\end{displaymath}
\item We  have that $\Delta_k g=0$.
\item The function $g(z)$ has
at most linear exponential growth at all the cusps of $\Gamma$.
\end{enumerate}

Theorem  \ref{transtheorem1} implies 
\begin{corollary} \label{corollarytransf2}
We have
\begin{multline*} 
R_2 \left(-\frac{1}{z};w \right) \,%
 e^{\frac{\pi i}{12  z}}  
= \frac{1}{(- i z)^{\frac{1}{2}}\cdot \eta(z)} 
\left( 
\frac{1}{1- e^{2 \pi i w} }
-  z \, \frac{e^{- 3 \pi i w^2 z - \pi i w + \pi i w z
}}{1 -  e^{2 \pi i  w  z }}  \right)   \\ 
+   i ( -i z)^{\frac{1}{2}}  e^{ - \frac{\pi i z}{12}  - 3 \pi i w^2 z - \pi i  w + \pi i w z}
R_2 \left(z; w z \right)
  - (-i z)^{-\frac{1}{2}}    \left( I_{1,0}^{+}\left(\frac{i}{z} ;w \right) +  I_{1,0}^{-}\left(\frac{i}{z} ;w \right)   \right).
\end{multline*}
\end{corollary}
Now  let  
\begin{eqnarray*}
J^{\pm}(z;w):= 
e^{ \pm \frac{\pi i}{6}}
\int_{\R} \frac{e^{ -\frac{3 \pi i x^2}{z} \pm \frac{\pi i x}{z}} }{ 1- 
e^{2 \pi i w \pm \frac{\pi i}{3} \pm \frac{2 \pi ix}{z}}  }
  dx.
\end{eqnarray*}
We show that $L\left( J^+(z;w)+J^-(z;w) \right)$ can be written as a certain Mordell type integral.
\begin{lemma} \label{lemmathetaint}
We have
\begin{eqnarray*}
L\left( J ^+(z;w)+J^-(z;w) \right) 
=  \frac{\sqrt{3}z^2}{2 \pi} \int_{0}^{\infty} 
\frac{\eta(iu)}{( - i (iu +z))^{\frac32}}\, du.
\end{eqnarray*}
\end{lemma}
\begin{proof}
Via analytic continuation it is enough to show the claim for   $z = it$.  The 
substitution $x \mapsto -  \frac{x}{t}$  gives
\begin{eqnarray} \label{integralcont}
J^{\pm} (it;w)  
= t \int_{\R} e^{- 3 \pi t x^2}  \frac{e^
{\pm \frac{\pi i}{6} \mp \pi x}}{1 - e^{2 \pi i w \pm \frac{\pi i }{3} \mp 2 \pi x }
}
dx.
\end{eqnarray}
We rewrite the integrand using the Mittag-Leffler theory of partial fraction decompositions.
\begin{eqnarray} \label{Mittag}
\sum_{\pm} \frac{e^ {\pm \frac{\pi i}{6} \mp \pi x}
 }{1 - e^{2 \pi i w \pm \frac{\pi i }{3} \mp 2 \pi x }  }
 = \mp
 \frac{ e^{ - \pi i w}  }{2 \pi}
 \sum_{m \in \Z}   
 (-1)^m \sum_{\pm} 
 \frac{1}{ x- i  \left( m + \frac{1}{6} \pm w \right) }.
\end{eqnarray} 
To see this identity, we first observe that both  
sides  of  the  identity are meromorphic functions of $x$ and it  is not hard to see that they have the same poles and residues.  
Moreover we will need later that 
the right-hand side is absolutely convergent. Thus in the following we are allowed to interchange summation, integration, and differentiation. 
Here it is important to consider the terms coming from the $\pm$ terms combined. 
Since both sides as a function of $x$ have period $i \Z$, it can be easily seen that the difference of both sides  is bounded on  $\C$ and  thus constant. Letting for fixed imaginary part  the real part of $x$ go to infinity
gives that this constant must be $0$.  
This gives 
\begin{eqnarray*}
J^{\pm}(it;w) = \mp
\frac{e^{- \pi i w} }{2 \pi}  \, t
\sum_{m }  (-1)^m 
\int_{\R} \frac{e^{-3 \pi t x^2}}{x- i \left(m + \frac{1}{6} \pm w  \right)}\, dx. 
\end{eqnarray*} 
It is not hard to see that
\begin{eqnarray*}
L \left( 
\frac{e^{ - \pi i w}}{ i w + x - i \left( m + \frac{1}{6}\right)} -
\frac{e^{ - \pi i w}}{- i w + x - i \left( m + \frac{1}{6}\right)}  \right)
=  - \frac{1}{\pi \left(x- i \left( m + \frac{1}{6}\right)  \right)^2}.
\end{eqnarray*}
Thus 
\begin{eqnarray}   \label{lidentity} 
L \left( J^+(it;w) +  J^-(it;w) \right) 
= 
 \frac{ t}{2 \pi^2}
\sum_{m}  (-1)^m \int_{\R} \frac{e^{ -3 \pi tx^2}}{\left(x- i\left( m+ \frac{1}{6}\right) \right)^2}\, dx.
\end{eqnarray}
For   $s \not= 0$ we have the identity
$$
\int_{\R} \frac{e^{- 2 \pi t x^2 }}{(x-is)^2} dx
= - \sqrt{2} \pi t \int_{0}^{\infty} 
\frac{e^{- 2 \pi u s^2}}{(u+t)^{\frac32} } du.
$$
Thus  
$$
L \left( J^+(it;w) +  J^-(it;w) \right) 
=
 \frac{\sqrt{3} (it)^2}
{2\pi}
\sum_{m} (-1)^m \int_{0}^{\infty} 
\frac{e^{- 3 \pi u  \left(m + \frac{1}{6} \right)^2 }}{(u+t)^{\frac32}} du,
$$
which gives the claim of the lemma.
\end{proof}  
Combining Corollary \ref{corollarytransf2}  with Lemma  \ref{lemmathetaint}  gives,
 as in the proof of Corollary, 
\ref{transcorollary1}
 \begin{corollary} \label{transcorollary3}
\begin{multline*} 
\mathcal{R} \left(-\frac{1}{z} \right) \,  
=   
\frac{2 \sqrt{6}}{\eta(z/24)}
\left( -
\frac{ (-i z)^{\frac{1}{2}}}{32 \pi}  
+ \frac{(-i z)^{\frac{3}{2}}}{24^3}\,   
+ \frac{1}{24(-i z)^{\frac{1}{2} }} \right)  +
\frac{(-i  z)^{\frac{3}{2}}}{48\sqrt{6}}  \mathcal{R} \left( \frac{z}{576}\right)  
\\ 
+ \frac{(-i z)^{\frac{3}{2}}}{16 \sqrt{3}\pi}
\int_{0}^{\infty} 
\frac{\eta(iu/24)}{(-i(iu+ z))^{\frac32}}\, du.
\end{multline*}
\end{corollary}
We next study the transformation law of  the non-holomorphic part  of $\mathcal{M}(z)$ and show that under inversion it introduces the same error integral as $\mathcal{R}(z)$.
\begin{lemma}  \label{lemmanonholtransf}
We have 
\begin{eqnarray*}
\mathcal{N}(z + 1)&=&  \mathcal{N}(z),\\
\mathcal{N} \left( - \frac{1}{z}\right) 
&=&  
 - \frac{  i( - i z)^{\frac32}}{16 \sqrt{3}\pi} 
  \int_{ - \bar{z }}^{i \infty}
\frac{\eta( \tau/24)}{(- i ( \tau + z))^{\frac32}}  d \tau
 +   \frac{(-i z)^{\frac32}}{16 \sqrt{3} \pi}
 \int_{0}^{\infty} 
 \frac{\eta\left( \frac{it}{24} \right)}{(-i (it + z) )^{\frac32}}\, dt
 .
\end{eqnarray*}
\end{lemma}
\begin{proof}
The first claim follows directly  since $\eta(24 \tau)$ is translation-invariant.
Moreover 
\begin{eqnarray*}
\mathcal{N}  \left( -\frac{1}{z} \right) 
= - 
 \frac{i}{4\sqrt{2} \pi}  \int_{\frac{1}{\bar z}}^{ i \infty}  
\frac{\eta(24\tau)}{\left( -i \left( -\frac{1}{z} +\tau \right)\right)^{\frac32}} d\tau.
\end{eqnarray*}
Making the change of variables $\tau \mapsto - \frac{1}{\tau}$ and using the transformation law of the  $\eta$-function, gives 
\begin{eqnarray*}
\mathcal{N} 
\left( -\frac{1}{z}  \right) 
= 
 - \frac{i(- i z)^{\frac32} }{ 16 \sqrt{3} \pi }\int_{ - \bar  z}^{ i \infty} \frac{\eta(\tau/24)}{(-i (\tau+z))^{\frac32}}\, d \tau
+  \frac{i(-i z)^{\frac32}}{ 16 \sqrt{3}\pi} \int_{0}^{ i \infty} 
\frac{\eta(it/24)}{(-i (z+it))^{\frac32}} \, dt,
\end{eqnarray*} 
as claimed.
\end{proof}
Combining Corollary \ref{transcorollary3}, Lemma \ref{lemmanonholtransf}, equation
(\ref{E2transf}), and  the transformation law of the $\eta-$function gives that $\mathcal{M}(z)$ transforms  
correctly under  $\Gamma_0(576)$ with Nebentypus character $\chi_{12}$.
To see that $\mathcal{M}(z)$ is annihilated by $\Delta_{\frac32}$,  we
write 
$
\Delta_{\frac32} = - 4 y^{\frac{1}{2} } 
\frac{\partial}{\partial z} y^{\frac{3}{2} } 
\frac{\partial}{\partial \bar z}.
$
Clearly 
$$
\frac{\partial}{\partial \overline z}
\left( 
\mathcal{R}(z) - \frac{1}{24 \eta(24 z)} 
+ \frac{E_2(24 z)}{8 \eta(24  z)}
\right)
=0.
$$
Moreover
$$
-\frac{\partial}{\partial \overline z} \mathcal{N}(z)
= \frac{i \eta(- 24 \overline z )}{16 \pi y^{\frac{3}{2} }} .
$$
Thus $y^{\frac{3}{2}} \frac{\partial}{\partial \overline z} \mathcal{N}(z)$
 is antiholomorphic, which implies that
 $\Delta_{\frac{3}{2}}\left( \mathcal{N}(z)\right)=0$. The claim about the 
 exponential growth in all cusps follows as in \cite{BO2}
 \section{Proof of Theorem \ref{maintheorem1}} \label{AsympSection}
Here we use the circle method and prove Theorem \ref{maintheorem1}.
First we estimate $I_{k,\nu}(z)$ and certain  Kloosterman sums.
\begin{lemma} \label{lemmaintest1}
Assume that $n \in \N,\, \nu \in \Z$, and $z:=\frac{k}{n}-k \Phi i$, $-\frac{1}{k(k+k_1)}\leq \Phi \leq \frac{1}{k(k+k_2)}$, where $\frac{h_1}{k_1}  < \frac{h}{k}  <\frac{h_2}{k_2} $ are adjacent Farey fractions in the Farey sequence of order $N$, with $N:=\lfloor n^{\frac{1}{2}}\rfloor$. Then
\begin{eqnarray*}
z^{\frac{1}{2}} \cdot  I_{k,\nu}(z)  \ll   
\frac{n^{\frac{1}{4}}}{\left\{\frac{\nu}{k} - \frac{1}{6k}  \right\}^2}
,
\end{eqnarray*}  
 where for $x \in \R$  we let $\{x \}:=x-\lfloor x \rfloor$. 
\end{lemma}
\begin{proof} 
We proceed similarly as in  \cite{An1} and \cite{Br}.
We write $\frac{\pi z}{k}=Ce^{iA}$. Then $|A|<\frac{\pi}{2}$ since $\re(z)>0$. 
Making the substitution $\tau=\frac{\pi z x}{k}$ 
gives 
\begin{eqnarray*}
z^{\frac{1}{2}} \cdot  I_{k,\nu}(z)= 
\frac{k}{\pi z^{\frac{1}{2}}}
\int_{S} 
e^{-\frac{3 k \tau^2}{\pi z}}  
\frac{\cosh \left( \frac{\pi i \nu}{k} - \frac{\pi i}{6k}   - \tau \right)}{\sinh^2 \left( \frac{\pi i \nu}{k}  - \frac{\pi i}{6k} - \tau \right)}\, d \tau,
\end{eqnarray*}
where  $\tau$ runs on the ray through $0$ of elements with argument $\pm A$.
One can easily  see that  for  $0 \leq t \leq A$
\begin{eqnarray*}
\left| e^{-\frac{3 k R^2 e^{2 i t}}{\pi z}} 
\frac{\cosh \left( \frac{\pi i \nu}{k} - \frac{\pi i}{6k}   \pm R\, e^{ it}  \right)}{\sinh^2 \left( \frac{\pi i \nu}{k}  - \frac{\pi i}{6k} \pm R\, e^{it} \right)}\,
\right| \to 0  \qquad \text{as } R \to \infty.
\end{eqnarray*}
Moreover the integrand  only has poles at points $ir$ with $r \in \R \setminus\{0 \}$. 
Shifting the path of integration
 to the real line  gives
\begin{eqnarray*}
z^{\frac{1}{2}}\cdot  I_{k,\nu}(z)= 
\frac{k}{\pi z^{\frac{1}{2}}}
\int_{\R} e^{-\frac{3 k t^2}{\pi z}}\cdot  
\frac{\cosh \left( \frac{\pi i \nu}{k} - \frac{\pi i}{6k}   - t \right)}{\sinh^2 \left( \frac{\pi i \nu}{k}  - \frac{\pi i}{6k} - t  \right)}
dt.
\end{eqnarray*}
We have 
\begin{eqnarray*}
\left|  e^{-\frac{3 k t^2}{\pi z}}\right| = e^{-\frac{3 k}{\pi} \re \left(\frac{1}{z} \right)t^2},\\
\left| \cosh \left(\frac{\pi i \nu}{k} -\frac{\pi i }{6k}- t \right)  \right|\leq e^t,
\end{eqnarray*}
\begin{displaymath}
\left| \sinh \left(\frac{\pi i \nu}{k} -\frac{\pi i }{6k}- t  \right)  \right|^2
\gg
\left\{
\begin{array}{ll}
e^{2t}&\text{if } t \geq 1,\\[1ex] 
\left\{ \frac{\nu}{k} - \frac{1}{6k} \right\}^2&\text{if } t \leq 1.
\end{array}
\right.
\end{displaymath}  
Thus 
\begin{eqnarray*}
z^{\frac{1}{2}} \cdot  I_{k,\nu}(z)
\ll \frac{k}{ 
  \left\{ \frac{\nu}{k} -\frac{1}{6k}\right\}^2    |z|^{\frac{1}{2}}  }    \int_{\R} e^{-\frac{3k}{\pi}  t^2 \re\left(\frac{1}{z} \right) } dt.
\end{eqnarray*}
Making the change of variables $t'= \sqrt{\frac{3 k \text{Re} \left( \frac{1}{z} \right)}{\pi} } \, t$ and   using  that 
$\re\left( \frac{1}{z}\right)^{-\frac{1}{2}} \cdot |z|^{-\frac{1}{2}} \ll n^{\frac{1}{4}} \cdot k^{-\frac{1}{2}}$ 
gives the claim of the lemma.
\end{proof}
We next give estimates for certain  sums of Kloosterman type \cite{An1, Br}. 
\begin{lemma} \label{lemmakloostest1}  
Let $n,m \in \Z$, $0 \leq \sigma_1 < \sigma_2 \leq k$, and  $D \in \Z$ with $(D,k)=1$. Then we have 
\begin{eqnarray} 
\label{kloost}
\sum_{h\pmod k^* \atop \sigma_1 \leq D h' \leq \sigma_2} 
\omega_{h,k}  \cdot e^{\frac{2 \pi i }{k} (hn+h'm )}  
\ll \gcd(24n+1,k)^{\frac{1}{2}}  \cdot
k^{\frac{1}{2}+\epsilon}.
\end{eqnarray}
\end{lemma}
\begin{proof}[Proof of Theorem \ref{maintheorem1}]
Using Lemma \ref{lemmaintest1}, Lemma \ref{lemmakloostest1}, and Corollary \ref{transcorollary1}, we now prove Theorem \ref{maintheorem1} 
 using the Hardy-Ramanujan method.
By Cauchy's Theorem we have  for $n>0$
\begin{eqnarray*}
\eta_2(n )
=\frac{1}{ 2 \pi i }  
\int_{C} \frac{R_2(q)}{q^{n+1}} \ dq,
\end{eqnarray*}
where $C$ is an arbitrary path inside the unit  circle surrounding $0$ counterclockwise.
Choosing for $C$ the circle with radius $e^{-\frac{2\pi}{n}}$ and as a parametrization $q= e^{-\frac{2 \pi}{n} + 2 \pi i t }$ with $0 \leq t \leq 1$,  gives
\begin{eqnarray*}
\eta_2(n)
= \int_{0}^{1} 
R_2 \left( e^{-\frac{2 \pi}{n} + 2 \pi i t}\right) \cdot e^{2 \pi- 2 \pi i n t}  \ dt. 
\end{eqnarray*} 
Define \begin{eqnarray*}
\vartheta_{h,k}' :=\frac{1}{k(k_1+k)},\quad 
\vartheta_{h,k}'' :=\frac{1}{k(k_2+k)}, 
\end{eqnarray*}
where $\frac{h_1}{k_1}  < \frac{h}{k}  <\frac{h_2}{k_2} $ are adjacent Farey fractions in the Farey sequence of order $N:=\left\lfloor n^{1/2} \right \rfloor$.  From the theory of Farey fractions it is known that 
\begin{eqnarray} \label{farey}
\frac{1}{k+k_j} \leq \frac{1}{N+1} \qquad (j =1,2).
\end{eqnarray}
We decompose the path of integration in paths along the Farey arcs 
$-\vartheta_{h,k}' \leq \Phi \leq \vartheta_{h,k}''$, where $\Phi:=t-\frac{h}{k}$.
This   gives
\begin{eqnarray*}
\eta_2 \left(n \right) = 
\sum_{h,k} 
e^{- \frac{2 \pi i hn}{k}}
\int_{-\vartheta_{h,k}'}^{\vartheta_{h,k}''}
R_2\left(e^{\frac{2 \pi i }{k}(h+iz) } \right)  \cdot 
e^{\frac{2 \pi n z}{k}} \ d\Phi,
\end{eqnarray*}
where $z:=\frac{k}{n}- k \Phi i$. Corollary \ref{transcorollary1} gives
\begin{multline*}
\eta_2(n) = 
\sum_{h,k} e^{-\frac{2 \pi i hn}{k}}\, \omega_{h,k}  \int_{- \vartheta_{h,k}'}^{\vartheta_{h,k}^{''}}
 e^{\frac{2 \pi nz}{k} }   e^{\frac{\pi}{12 k} (z^{-1} -z)}
 \, z^{\frac{1}{2}}  
 \left(  
 \frac{1}{(q_1;q_1)_{\infty}}
 \left( 
-\frac{3k}{4\pi z} + \frac{1}{24z^2} + \frac{1}{24}\right) \right.  \\ 
\left.
+ \frac{1}{z^2} R_2(q_1)  \right)  d\Phi   
-  \sum_{h,k} \frac{\omega_{h,k}}{k} e^{ - \frac{2 \pi i hn}{k} }
\sum_{\nu \pmod k} 
(-1)^{\nu} e^{\frac{\pi i h'}{k} (-3 \nu^2+ \nu) }  
 \int_{- \vartheta_{h,k}'}^{\vartheta_{h,k}^{''}}   
 z^{\frac{1}{2}} \,
 I_{k,\nu}(z)
e^{\frac{2 \pi (n-1/24)z}{k}}\, d \Phi,
 \end{multline*}
 and we denote the  two  summands by $\sum_1$ and $\sum_2$.
 We first treat $\sum_1$ and start with   the contribution coming from 
 $$
 \frac{1}{24z^2 \, (q_1;q_1)_{\infty}} 
 =: \frac{1}{z^2} \left(
 \frac{1}{24} + \sum_{r>0} a(r)\, q_1^r \right).
 $$
 We consider  the constant term and the term arising from $r \geq 1$ separately  
  since they contribute to the main term and to the  error term, respectively.  We denote the associated sums by $S_1$ and $S_2$, respectively and first estimate $S_2$. Throughout we need the easily verified facts  that  $\re(z)=\frac{k}{n}$, $\re\left( \frac{1}{z}\right)> \frac{k}{2}$,
$|z|  \geq \frac{k}{n}$, and $\vartheta_{h,k}'+\vartheta_{h,k}^{''} \leq
\frac{2}{k(N+1)}$. 
Since $k_1, k_2 \leq N$, we can  write 
\begin{eqnarray} \label{pathint}
\int_{- \vartheta_{h,k}'}^{\vartheta_{h,k}^{''}}
=
\int_{-\frac{1}{k(N+k)}}^{\frac{1}{k(N+k)}}  
+ \int_{-\frac{1}{k(k_1+k)}}^{-\frac{1}{k(N+k)}} 
+ \int_{\frac{1}{k(N+k)}}^{\frac{1}{k(k_2+k)}} 
\end{eqnarray}
and we denote the associated sums by $S_{21}$, $S_{22}$, and $S_{23}$, respectively.    
We first consider $S_{21}$.  Lemma \ref{lemmakloostest1} gives
\begin{multline*}
S_{21} \ll
\left| \sum_{r=1}^{\infty} 
 a(r)
 \sum_{k}   \sum_{h}  \omega_{h,k} 
 \cdot 
e^{ - \frac{2 \pi i hn}{k} + \frac{2 \pi i r h' }{k}} 
\int_{- \frac{1}{k(N+k)}}^{   \frac{1}{k(N+k)}}
z^{-\frac{3}{2}} \cdot 
e^{-\frac{2 \pi}{kz}\left(r-\frac{1}{24}\right)+\frac{2 \pi z}{k}\left(n-\frac{1}{24} \right)}  \ d \Phi  \right| \\ 
\ll  
\sum_r |a(r) | e^{ - \pi r }
n  \sum_{k} (24n-1,k)^{\frac{1}{2}} \, k^{-2+ \epsilon} 
\ll n.
\end{multline*}
Since  $S_{22}$ and $S_{23}$ are estimated similarly   we only consider $S_{22}$. 
Writing 
\begin{eqnarray*}
\int_{-\frac{1}{k(k+k_1)}}^{-\frac{1}{k(N+k)}}
= \sum_{l=k_1+k}^{N+k-1} \int_{-\frac{1}{kl}}^{-\frac{1}{k(l+1)}},
\end{eqnarray*}
we see 
\begin{multline} \label{est4}
S_{22} \ll 
 \left|   \sum_{r=1 }^{\infty} a(r)
  \sum_{k} 
\sum_{l=N+1}^{N+k-1} \int_{-\frac{1}{kl}}^{-\frac{1}{k(l+1)}} 
z^{-\frac{3}{2}} \cdot 
e^{-\frac{2 \pi}{kz}\left(r-\frac{1}{24}\right)+\frac{2 \pi z}{k}\left(n-\frac{1}{24} \right)} 
d\Phi 
\sum_{h\atop N < k+k_1 \leq l}  \omega_{h,k}   \cdot e^{- \frac{2 \pi i hn}{k} + \frac{2 \pi i r h'}{k}} \right|. 
\end{multline}
It follows from the theory of Farey fractions that
\begin{eqnarray*}
k_1 \equiv - h' \pmod k, \quad k_2 \equiv h' \pmod k,  
\\     
N-k<k_1 \leq N, \quad N-k < k_2 \leq N.
\end{eqnarray*}
Thus (\ref{est4}) can be estimated similarly as $S_{21}$ 
  using Lemma \ref{lemmakloostest1}.   
  
In the same way we decompose the terms   in 
$ 
-
\frac{3k}{4 \pi z(q_1;q_1)_{\infty} },\, \frac{1}{24(q_1;q_1)_{\infty}},
$
and $\frac{R_2(q_1)}{z^2}$.  
One can show that the term with a  positive exponent in the Fourier expansions also introduces an error of  $O\left( n^{ 1+ \epsilon} \right)$. 
Thus 
\begin{eqnarray*} 
\sum_1=\sum_{h,k} e^{-\frac{2 \pi i hn}{k}}\, \omega_{h,k}  \int_{- \vartheta_{h,k}'}^{\vartheta_{h,k}^{''}}
 e^{\frac{2 \pi \left(n -\frac{1}{24} \right)z}{k} +\frac{ \pi}{12kz} }    
 \left( 
-\frac{3k}{4\pi z^{\frac{1}{2}}} + \frac{1}{24z^{\frac{3}{2}}} + \frac{z^{\frac{1}{2}}}{24} 
\right)   
+ O\left(n^{1+ \epsilon}  \right). 
 \end{eqnarray*}
 We next write the path of integration in a symmetrized way
 
  \begin{eqnarray*}
\int_{- \vartheta_{h,k}'}^{\vartheta_{h,k}^{''}}
=
\int_{-\frac{1}{kN}}^{\frac{1}{kN}}  -
 \int_{-\frac{1}{kN}}^{-\frac{1}{k(k+k_1)}} 
- \int_{\frac{1}{k(k+k_2)}}^{\frac{1}{kN}}
\end{eqnarray*}
and denote the associated sums by $S_{11}$, $S_{12}$, and $S_{13}$, respectively.
The sums  $S_{12}$ and $S_{13}$ contribute to the error term. 
Since they  have a similar  shape, we only consider $S_{12}$.
We again only estimate the error arising  from  the second  summand since it can be shown the other terms lead to an error of at most that size. 
Writing
\begin{eqnarray*}
\int_{-\frac{1}{kN}}^{-\frac{1}{k(k_1+k_1)}}
= \sum_{l=N}^{k+k_1-1} 
\int_{-\frac{1}{kl}}^{-\frac{1}{k(l+1)}}
\end{eqnarray*}
gives
\begin{eqnarray} \label{S12}
S_{12} \ll  \left| 
 \sum_{k} 
\sum_{l=N}^{N+k-1} 
\int_{-\frac{1}{kl}}^{-\frac{1}{k(l+1)}} 
z^{-\frac{3}{2}}
e^{\frac{ \pi}{12kz}+\frac{2 \pi z}{k}\left(n-\frac{1}{24} \right)} 
d\Phi  
\sum_{h\atop l < k+k_1-1 \leq N+k-1}  \omega_{h,k}     \cdot e^{- \frac{2 \pi i hn}{k} } \right|
. 
\end{eqnarray} 
Using that $\re(z)=\frac{k}{n}$, $\re\left( \frac{1}{z}\right)<k$, and $|z|^2 \geq \frac{k^2}{n^2}$, (\ref{S12}) can  by Lemma \ref{lemmakloostest1} be estimated as before   against $O\left(n^{1+ \epsilon}   \right)$.
Thus  
\begin{eqnarray} \label{sum1est}
\sum_1=  \sum_{k} 
A_{k}(n)
\int_{- \frac{1}{kN}}^{\frac{1}{kN}} 
\left( 
- \frac{3k}{4\pi z^{\frac{1}{2}}} +  \frac{1}{24z^{\frac{3}{2}}} 
+ \frac{z^{\frac{1}{2}}}{24}
\right)  
e^{\frac{2 \pi z}{k}\left(n-\frac{1}{24}  \right)+\frac{\pi}{12k z }}  \
d \Phi
+ O\left( n^{1+ \epsilon} \right).
\end{eqnarray}
To finish the estimation of $\sum_1$,  we  consider integrals of the form
 \begin{eqnarray*}
I_{k,r}:= \int_{-\frac{1}{kN}}^{\frac{1}{kN}} 
z^{r} \cdot 
e^{\frac{2 \pi}{k}\left(z \left(n - \frac{1}{24} \right) +\frac{1}{24z}\right)} d\Phi
,
\end{eqnarray*}  
where $r \in \{- \frac{3}{2},-\frac{1}{2}, \frac{1}{2} \}$. 
Substituting $z=\frac{k}{n}- i k \Phi$   gives 
\begin{eqnarray} \label{int}
I_{k,r}= \frac{1}{ki} 
\int_{\frac{k}{n}-\frac{i}{N}}^{\frac{k}{n}+\frac{i}{N}} 
z^{r} \cdot 
e^{\frac{2 \pi}{k}\left(z \left(n - \frac{1}{24} \right) +\frac{1}{24z}\right)} \ dz 
.
\end{eqnarray}
We   denote the circle through $\frac{k}{n} \pm \frac{i}{N}$ and tangent to the imaginary axis at $0$ by $\Gamma$. 
If $z=x+iy$, then $\Gamma$ is given by $x^2 + y^2 = \alpha x$, with $\alpha = \frac{k}{n} +\frac{n}{N^2k}$.
Using the fact that $2 > \alpha >\frac{1}{k}$, $\re(z)\leq \frac{k}{n}$, and $\re\left( \frac{1}{z}\right)<k$ on the smaller arc we can show that the integral along the smaller arc is in $O \left(n^{-\frac{1}{4}} \right)$. Moreover the path of integration in (\ref{int}) can be changed by Cauchy's Theorem into the larger arc of $\Gamma$.  
Thus 
\begin{eqnarray*}
I_{k,r}= \frac{1}{ki} 
\int_{\Gamma} 
z^{r} \cdot
e^{\frac{2 \pi}{k}\left(z \left(n - \frac{1}{24} \right) +\frac{1}{24z}\right)}\ dz 
+ O \left( n^{-\frac{1}{4}}\right).
\end{eqnarray*}    
Making the substitution $t = \frac{\pi}{12kz}$ gives 
\begin{displaymath}
I_{k,r} 
=
\frac{2 \pi}{k} \left( \frac{\pi}{12k} \right)^{1+r}
\frac{1}{2 \pi i } 
\int_{\gamma - i \infty}^{\gamma+ i \infty}
t^{-2-r} \cdot e^{t+\frac{\alpha}{t}} \ dt + O\left(n^{-\frac{1}{4}} \right),
\end{displaymath}
where $\gamma \in \R^+$ and $\alpha :=\frac{ \pi^2 }{144k^2}(24n-1)$. 
Now
$$
\frac{1}{2 \pi i} \int_{\gamma - i \infty }^{\gamma+ i \infty}
t^{-(2+r) } \, e^{ t + \frac{\alpha}{t}}\, dt 
= \alpha^{- \frac{r+1}{2} }\, I_{r+1 } \left(2 \sqrt{\alpha}\right).
$$
Thus   
$$
I_{k,r} =   
\frac{2\pi}{k}
 (24n-1)^{ - \frac{r+1}{2}} I_{r+1 }\left(\frac{\pi \sqrt{24n-1}}{6k} \right) 
+ O \left( n^{-\frac{1}{4} + \epsilon} \right).
$$  
Plugging this into (\ref{sum1est}) gives 
\begin{multline*}
\sum_1
= \sum_k A_k(n)
\left(     - \frac{3}{2
  (24n-1)^{\frac{1}{4}}}
 I_{\frac{1}{2}} \left( \frac{\pi}{6k} \sqrt{ 24n-1} \right) 
+ \frac{\pi(24n-1)^{\frac{1}{4}} }{12k} \right. \\ 
\left. \, I_{-\frac{1}{2}} \left(\frac{\pi}{6k} \sqrt{24n-1} \right)
+  \frac{\pi}{12k (24n-1)^{\frac{3}{4}}} \, I_{\frac{3}{2}} \left(\frac{\pi}{6k} \sqrt{24n-1} \right)
\right) 
+ O \left( n^{1+ \epsilon }\right).
\end{multline*}
We next consider $\sum_2$.  Using  Lemmas  \ref{lemmaintest1} and  \ref{lemmakloostest1},  we obtain
\begin{eqnarray*}
\sum_{2} \ll 
n^{-\frac{1}{4}} \sum_k k^{\frac{3}{2} + \epsilon}
(24n-1,k)^{\frac{1}{2}}
 \sum_{ \nu \pmod k} 
\frac{1}{\left(k \left\{\frac{\nu}{k} - \frac{1}{6}\right\} \right)^2}
\ll n^{-\frac{1}{4}} \sum_k k^{\frac{3}{2} + \epsilon } 
\sum_{\nu =1 }^{k} 
\frac{1}{\nu^2} 
\ll n^{\frac12 + \epsilon}.
\end{eqnarray*}  
Combining the estimates for $\sum_1$ and $\sum_2$  gives Theorem \ref{maintheorem1}.
\end{proof}

Corollary \ref{corasymp1} can  be concluded using that  for $x \to \infty$
$$
I_{\alpha} (x) \sim \frac{1}{\sqrt{2\pi x}} \, e^x.
$$  
\section{Proof of Theorem \ref{congtheorem}} \label{CongSection}
Here we prove Theorem \ref{congtheorem}.
First we need to know on which arithmetic progressions the non-holomorphic part of $\mathcal{M}(z)$ is supported.
Similarly as in \cite{BO2} one can show.
\begin{lemma}
We have
\begin{equation*}
\mathcal{N}(z)
=   \frac{1}{4 \sqrt{\pi}} \sum_{k \in \Z}  (-1)^k 
(6k+1) \Gamma  \left(-\frac12; 4\pi (6k+1)^2 y \right)
q^{- (6k+1)^2 },
\end{equation*}
where 
$\Gamma(\alpha;x):= \int_{x}^{\infty} e^{-t}\, t^{\alpha-1}\, dt$ 
is the incomplete Gamma-function.
\end{lemma}
We next recall results from \cite{BO2}.  For  this  let $0 <a<c$ and define
\begin{eqnarray*}
D\left( \frac{a}{c};z\right) &:=& - S \left( \frac{a}{c};z\right) +q^{-\frac{\ell_c}{24} }\, R \left(\zeta_c^a;q^{ \ell_c} \right),\\
S\left( \frac{a}{c};z\right) &:=& - \frac{i \sin \left( \frac{\pi a}{c}\right) \ell_c^{\frac{1}{2} }}{\sqrt{3}}\,
\int_{-\bar z}^{i \infty}  
\frac{\Theta\left( \frac{a}{c};\ell_c \tau\right)}{ \sqrt{ i (\tau+z) } } d \tau,
\end{eqnarray*}
where  $\ell_c:= \lcm (2c^2,24)$, $\zeta_c:=e^{\frac{2 \pi i}{c}}$, 
 and  $\Theta\left(\frac{a}{c};\tau \right)$ is a certain weight $\frac{3}{2}$ cuspidal theta function (for the definition    see \cite{BO2}). 
It turns out (see Theorems 1.1 and 1.2 of  \cite{BO2}) that $D\left(\frac{a}{c};z \right)$  is a harmonic 
weak Maass form of weight $\frac{1}{2}$.  

Next we observe, using orthogonality of roots of unity, that 
\begin{eqnarray*}
\sum_{n=0}^{\infty}
NF_2(r,t;n)\, q^n 
= \frac{1}{t} R_2(q) + \frac{1}{t}\sum_{j=1}^{t-1}
\zeta_t^{-rj}
R_2\left(\zeta_t^{j};q \right).
\end{eqnarray*}
Thus 
\begin{multline*}
\sum_{n=0}^{\infty}
NF_2(r,t;n)q^{ \ell_t-\frac{\ell_t}{24}  } 
= \frac{1}{t} \left(  \mathcal{M} \left( \frac{\ell_t z}{24} \right)  + \mathcal{N} \left( \frac{\ell_t z}{24} \right) 
+\frac{1}{24 \eta( \ell_t z )}
+ \frac{ E_2(\ell_t z)}{8 \eta(\ell_t )} \right) 
\\
+
\frac{1}{t}
\sum_{j=1}^{t-1} 
\frac{ \zeta_t^{-rj + 2j }}{ \left(1-\zeta_t^j \right)\left(\zeta_t^{3j} -1 \right)}
\left( 
D \left(\frac{j}{t};z \right) 
- D \left(\frac{2j}{t};z \right)
+ \left(S\left(\frac{j}{t};z \right)-S \left(\frac{2j}{t};z \right) \right)
\right).
\end{multline*}
Now fix a prime $p>3$, and let 
$\mathcal{S}_p:=\left\{ n \in \Z; \leg{24 \ell_tn -\ell_t }{p} = - \leg{\ell_t}{p} \right\}$.
From \cite{BO2}, one can conclude that the restriction of 
$
D \left(\frac{js}{t};z \right) 
$ $(s \in \{1,2\})$ to those coefficients lying in $\mathcal{S}_p$ is a weakly holomorphic modular form of weight $\frac12$ on $\Gamma_1(6 f_t^2 \ell_tp^4)$.
In a similar way we see that the restriction of   $\mathcal{M} \left( \frac{\ell_t z}{24} \right)$
to those coefficients in $\mathcal{S}_p$ is a weakly holomorphic modular form of weight $\frac32$ on $\Gamma_1\left(576 \ell_tp^4\right)$.
Moreover the restriction of  
$\frac{1}{24 \eta( \ell_t z)}$ 
to those coefficients in $\mathcal{S}_p$ is a weakly holomorphic modular form of weight $-\frac12$ on $\Gamma_1(24 \ell_tp^4)$.
Moreover by work  of Serre it is known that $E_2(z)$ is a $p$-adic modular form, e.g. we have  
\begin{eqnarray*}
E_2(z) 
\equiv E_{p+1}(z) \pmod{p}.
\end{eqnarray*}
From this one can conclude that    the restriction of 
$ \frac{ E_2(\ell_tz)}{ \eta(\ell_t z)} $ to $\mathcal{S}_p$
 is congruent to a weakly holomorphic modular form on    $\Gamma_1(24 \ell_tp^4)$ modulo $Q^j$.
Now the result can be concluded as in \cite{BO2}.  
 \section{A transformation law for $R_2^o(q)$} \label{TransoddSection}
 Throughout we assume the notation of Section \ref{TRANS}.
Define  for $w \in \C$ with Re$(w)$ sufficiently small
\begin{eqnarray*}
\mathcal{R}^o(q;w)= 
\mathcal{R}^o(z;w)&:=& \frac{1}{(q^2;q^2)_{\infty}}
\sum_{n \in \Z} 
\frac{(-1)^n\, q^{3n^2+n}}{1-e^{2 \pi i w}q^{2n+1}}
\end{eqnarray*}
and for a  function $g$ that is differentiable in a neighborhood of $0$,  let
\begin{eqnarray*}
L^o (g(w))&:=& - \frac{1}{8 \pi^2}\left[ \frac{\partial }{\partial w} e^{- 2 \pi i w } \frac{\partial}{\partial w} g
\right]_{w=0} .
\end{eqnarray*}
Then 
\begin{eqnarray*} 
L^o \left( \mathcal{R}^o(q;w)\right)=R^o(q).
\end{eqnarray*}
We show a transformation law for $\mathcal{R}^o(q;w)$ and then conclude a transformation law for $R^o(q)$ by applying $L^o$.
This transformation law turns out to be more complicated than the transformation law for $R(q;w)$.
We distinguish the cases $k$ odd and $k$ even.
Let 
$$
H^o_w(x):= 
\frac{1}{2 \sinh \left( - \pi i w - \frac{x}{2} \right)}.
$$
Clearly 
$$
H_w^o(x+ 2 \pi  i) = - H_w^o(x).
$$
Let 
\begin{eqnarray*}
h^o(z;w)=h^o(q;w)&:= &
\frac{1}{\left(q; q \right)_{\infty}}
\sum_{m \in \Z} 
(-1)^m 
\frac{q^{\frac{1}{2}(3m^2+m ) }}{1+ e^{\pi i w} q^m},\\
I_{ k,\nu}^{o,\pm} (z;w)&:= &
\int_{\R} e^{ - \frac{3 \pi zx^2}{2k} + \frac{3 \pi zx}{k} }
\, H_{w}^o \left( \pm
\frac{2 \pi i \nu}{k} - \frac{2 \pi zx}{k} \mp \frac{\pi i}{3k} 
\right)\, dx.
\end{eqnarray*}
\begin{theorem} \label{transtheoremodd1}
Assume that $k$ is odd and that $4|h'$. Then we   have 
\begin{multline*} 
\mathcal{R}^o(z;w)
= -
\frac{1}{\sqrt{2}} z^{- \frac{1}{2} } \omega_{2h,k}
 e^{\frac{\pi}{24kz} -\frac{2 \pi z}{3k} + \frac{\pi i h}{2k} }
e^{\frac{\pi i}{2} (k+ 1 + 3hk) }\,   
e^{2 \pi iw + \frac{3 \pi k w^2}{2z} - \frac{\pi w}{2z}  } h^o\left(q_1^{\frac{1}{2}}; \frac{i w}{z}\right)
\\   - 
\frac{1}{\sqrt{2}\, k} e^{\frac{\pi i}{2} k(3h+1) }\, \omega_{2h,k  } z^{ \frac{1}{2}} 
e^{ - \pi iw }  e^{\frac{\pi i h}{2k} -\frac{2 \pi z}{3k}}
\sum_{\substack{\nu \pmod k \\  \pm }} 
\pm  e^{ \frac{\pi i}{2k}\left( h'\left( - 3 \nu^2 + \nu \right) \mp  6 \nu \pm 1\right)}\, I_{k,\nu}^{o,\pm}(z;w).
\end{multline*}  
\end{theorem}
\begin{proof} 
We proceed as in the proof of Theorem \ref{transtheorem1}.
Via analytic continuation it is enough to show the claim for $w \in \C$ with Re$(w)\not =0$ sufficiently small and $z>0$ real. 
Define
$$  
\widetilde{R}^o(z;w)  :=  
\sum_{n \in \Z} 
\frac{(-1)^n\, q^{3n^2+n}}{1-e^{2 \pi i w}q^{2n+1}}.
$$
Then 
\begin{eqnarray} \label{tildefunction}
  \widetilde{R}^o(z;w) =
e^{- \pi i w} \, q^{- \frac{1}{2}} 
\sum_{n \in \Z} (-1)^n e^{\frac{6 \pi i (h+iz)n^2}{k} }
\, H_w^o \left( \frac{2 \pi i}{k}(h+iz)(2n+1) \right).
\end{eqnarray}
We   write $n=\nu + km$ with $\nu$ running modulo $k$ and $m \in \Z$.
This gives that   
\begin{multline} \label{poisson3}
\widetilde{R}^o(z;w)
=
e^{- \pi i w} \, q^{- \frac{1}{2}} 
\sum_{\nu \pmod k} 
(-1)^{\nu} \, 
e^{\frac{6 \pi i h \nu^2}{k} }
\sum_{m \in \Z} (-1)^m e^{-\frac{6 \pi z(\nu+ km)^2}{k} }\\
\, H_w^o \left(   \frac{2 \pi i h}{k}(2 \nu+1) - \frac{2 \pi z}{k}  \left(  2 \left( \nu +km \right) +1  \right)     \right).
\end{multline}
Applying 
Poisson summation to the inner sum  and making the change of variables $x \mapsto \nu +kx$ gives that the inner sum equals
$$
\frac{1}{k}
\sum_{n \in \Z}  \int_{\R} 
e^{-\frac{6 \pi z x^2}{k} }\\
\, H_w^o \left(   \frac{2 \pi i h}{k}(2 \nu+1) - \frac{2 \pi z}{k}  (  2x+1  )     \right) 
e^{\frac{\pi i}{k}(2n+1)(x - \nu)}
dx.
$$
In (\ref{poisson3}), we change for $n \leq -1$, $n \mapsto -n-1$, $x \mapsto -x$, and $\nu \mapsto - \nu$.
This gives 
\begin{multline} \label{poisson}
 \widetilde{R}^o(z;w)   =
e^{- \pi i w} \, q^{- \frac{1}{2}} 
 \frac{1}{k}
\sum_{\substack{\nu \pmod k\\ \pm}}  
(-1)^{\nu} \, 
e^{\frac{6 \pi i h \nu^2}{k} }
     \sum_{n \in \N_0}   \int_{\R} 
     e^{-\frac{6 \pi z x^2}{k} }\\
\, H_w^o \left(   \frac{2 \pi i h}{k}(\pm 2 \nu+1) - \frac{2 \pi z}{k}  ( \pm 2x+1  )     \right) 
e^{\frac{\pi i}{k}(2n+1)(x - \nu)}
dx.     
\end{multline} 
Making the change of variables $x \mapsto 2 x\pm 1$ and $\nu \mapsto \frac{k+1}{2}(\nu \mp 1)$, 
a lengthy calculation gives
\begin{multline*} 
 \widetilde{R}^o(z;w)
=
e^{- \pi i w} \, q^{\frac{1}{4}} \frac{1}{2k}
\sum_{\substack{\nu \pmod k \\ \pm}}  
e^{- \frac{\pi i}{2}(\nu \mp 1) + \frac{3 \pi i h}{2k} (\nu^2\mp 2 \nu) + \frac{\pi i}{2} ( k+1 + 3kh)(1- \nu^2)  }  
\\
\sum_{ n \in \N_0}  
(-1)^{n(\nu +1 ) } 
\int_{\R } 
e^{  - \frac{3 \pi z}{2k} (x^2 \mp 2x)}
e^{ \frac{\pi i}{2k} (2n+1)( x - \nu)}
H_w^o \left( 
\pm \frac{2 \pi ih \nu}{k} 
\mp \frac{2 \pi z x}{k}
\right)\, 
dx
\,   .     
\end{multline*}
We use the residue theorem to shift the path   of  integration through 
$$
\omega_n:= \frac{i}{6z} (2n+1).
$$
Using again the   function $S_w^{\pm}$ defined in (\ref{Sw}), one sees that poles of the integrand can only come from poles of 
$\sinh(\pi zx \mp \pi i k w)$,
thus lie in 
\begin{equation} \label{xm2}
x_m^{\pm} := \frac{i}{z} (m \pm kw).
\end{equation}
Moreover for fixed $m$, there is at most for one $\nu$ modulo $k$ a non-trivial  residue and this $\nu$ can be chosen as 
\begin{equation}  \label{nm2}
\nu_m:= -h' m.
\end{equation}
Thus poles may at most occur for $x_m$ and $\nu_m$ as in (\ref{xm2}) and (\ref{nm2}).
If we now shift the path of integration through $\omega_n$, we have to take those $x_m^{\pm}$ into account for which $n \geq 3m>0$, and for $x_m^+$ we consider additionally $m=0$.  
Denoting by $\lambda_{n,m}^{\pm}$ the residue of the integrand we get 
\begin{eqnarray*}
 \widetilde{R}^o(z;w)
= \sum_1 + \sum_2,
\end{eqnarray*}
\begin{eqnarray*}
\sum_1:= \frac{\pi i}{k} e^{ - \pi i w} \, q^{\frac{1}{4}}\,  \left(
\sum_{m \geq 0} \sum_{n \geq 3m} \lambda_{n,m}^+  
+ \sum_{m>0} \sum_{n \geq 3m} \lambda_{n,m}^{-}
\right),
\end{eqnarray*}
\begin{multline*}  
\sum_2
:=
e^{- \pi i w} \, q^{\frac{1}{4}} \frac{1}{2k}
\sum_{\substack{\nu \pmod k \\ \pm}}  
e^{- \frac{\pi i}{2}(\nu \mp 1) + \frac{3 \pi i h}{2k} (\nu^2\mp 2 \nu) + \frac{\pi i}{2} ( k+1 + 3kh)(1- \nu^2)  }  
\\
\sum_{ n \in \N_0}  
(-1)^{n(\nu +1 ) } 
\int_{\R + \omega_n } 
e^{  - \frac{3 \pi z}{2k} (x^2 \mp 2x)}
e^{ \frac{\pi i}{2k} (2n+1)( x - \nu)}
H_w^o \left( 
\pm \frac{2 \pi ih \nu}{k} 
\mp \frac{2 \pi z x}{k}
\right)\, 
dx
\,   .     
\end{multline*}
We first consider $\sum_1$. Since $h'$ is even, we have 
\begin{eqnarray*}
\lambda_{n,m}^{\pm} 
= \pm (-1)^n \frac{k\, e^{\frac{\pi i}{2}(k+1 + 3hk) - 
\frac{\pi i}{2}( \nu_m \mp 1) 
+ \frac{3 \pi i h}{2k}(\nu_m^2  \mp 2 \nu_m  ) - \frac{3 \pi z }{2k} \left( \left(x_{m}^{ \pm  }  \right)^2 \mp 2 x_{m}^{ \pm} \right) + \frac{\pi i}{2k} (2n+1)  ( x_m^{ \pm
} - \nu_m) }}{ 2 \pi z\,  \cosh ( \mp \pi i w - \frac{\pi i  h \nu_m}{k}  + \frac{\pi z x_m^{\pm}}{k})}.
\end{eqnarray*}
From this we  see that 
$$
\lambda_{n+1,m}^{\pm}
= - e^{\frac{\pi i}{k} ( x_m^{ \pm} - \nu_m)} \lambda_{n,m}^{\pm}
.
$$
Thus 
\begin{eqnarray*}
\sum_1  = \frac{\pi i}{k} e^{ - \pi i w} 
q^{ \frac{1}{4}}   
\left( 
\sum_{m \geq 0} 
\frac{\lambda_{3m,m}^+}{1 + e^{ \frac{\pi i}{k}  ( x_m^+  - \nu_m)}}
+ 
\sum_{ m > 0} 
\frac{\lambda_{3m,m}^-}{1 + e^{ \frac{\pi i}{k}  ( x_m^-  - \nu_m)}}
\right).
\end{eqnarray*}
A lengthy calculation calculation gives that this equals  
\begin{eqnarray*} 
- 
e^{2 \pi iw + \frac{3 \pi k w^2 }{2z} - \frac{\pi w}{ 2z } } 
q^{\frac{1}{4}} \, e^{\frac{\pi i}{2} (k+1+ 3hk)} 
\frac{1}{2z}    \left(q_1^{1/2};q_1^{1/2} \right)_{\infty}  
h^o \left( q_1^{1/2};\frac{iw}{z}\right).
\end{eqnarray*}
We next turn to $\sum_2$. As before we can see   that 
\begin{multline*}
\sum_2 = 
e^{- \pi i  w }\, q^{-\frac{1}{2}} 
\frac{1}{k} 
\sum_{\nu \pmod k} (-1)^{\nu} 
\, e^{\frac{6 \pi i h \nu^2}{k}} 
\sum_{n \in \Z}
\int_{\R + \frac{\omega_n}{2} }
e^{ - \frac{6 \pi z x^2 }{k}} \\
H_w^o \left( 
\frac{2 \pi i h}{k} (2 \nu+1) -\frac{2 \pi z}{k}(2x+1) 
\right)\, e^{ \frac{\pi i}{k} (2n+1)(x - \nu)}\, dx.
\end{multline*}
Substituting $x \mapsto x + \frac{\omega_n}{2}$ 
and writing  $n=3p+ \delta$ with $p \in \Z$ and
$\delta \in \{ 0,  \pm 1\}$ gives  
\begin{multline*}
\sum_2 = 
e^{- \pi  i  w }\, q^{-\frac{1}{2}} 
\frac{1}{k}  
\sum_{\nu \pmod k} (-1)^{\nu} 
\, e^{\frac{6 \pi i h \nu^2}{k}} 
\sum_{\delta, p}
e^{- \frac{\pi i}{k} (6p + 2 \delta +1) \nu - \frac{\pi}{24 kz}  (6p + 2 \delta +1)^2} \\
\int_{\R}
e^{ - \frac{6 \pi z x^2 }{k}}  
H_w^o \left( 
\frac{2 \pi i h}{k} (2 \nu+1) -\frac{2 \pi z}{k}(2x+1)   - \frac{\pi i}{3k} (6 p + 2 \delta +1)
\right)\,  dx.
\end{multline*}
We next change $\nu$ into $\frac{1+k}{2} (  - h' ( \nu+ p  ) -1)$ and $x \mapsto \frac{x-1}{2}$.
A lengthy calculation  gives
\begin{multline*}
\sum_2 =   - 
e^{ - \pi i w  } q^{ \frac{1}{4} } \frac{1}{2k}   
e^{    \frac{\pi i}{2}  k(1+3h) } 
\sum_{p,\delta, \nu} (-1)^{ p+\delta} 
q_1^{ \frac{1}{48} (6p+ 2 \delta+1)^2} 
e^{ \frac{\pi i}{2k} ((2\delta+1)- 6 \nu) }\\
e^{ \frac{\pi i h'}{2k} \left(- 3 \nu^2 + \nu(2 \delta+1)- \frac{1}{12}\left(2 \delta+1 \right)^2\right) }   
\int_{\R} 
e^{ - \frac{3 \pi z x^2}{2k}+ \frac{3 \pi zx}{k}} \, 
H_w^o
\left( 
\frac{2 \pi i \nu}{k} - \frac{2 \pi  zx }{k} - \frac{\pi i}{3k} ( 2 \delta+1)
\right)\, dx.
\end{multline*}
Now the integral is independent of $p$.
Moreover for $\delta=1$ the sum over $p$ vanishes since the $p$th and the $-(p+1)$th term cancel. 
For $\delta=0, - 1$ it equals 
$\eta \left( \frac{\pi i}{k} \left(h' + \frac{i}{z} \right)   \right)$.
Changing for $\delta=-1$, $\nu \mapsto -\nu$ gives
\begin{multline*} 
\sum_2 =   -
e^{ - \pi i w  } q^{ \frac{1}{4} } \frac{1}{2k}
  e^{\frac{\pi i}{2} k(3h+1)   -  \frac{\pi ih'}{24 k} } 
\eta \left( \frac{\pi i}{k} \left(h' + \frac{i}{z} \right)   \right) 
\sum_{\substack{\nu  \pmod k \\ \pm  } } 
\pm
e^{  \frac{\pi i}{2 k} \left(h' \left( - 3 \nu^2 + \nu \right)   \mp 6 \nu \pm 1    \right) }
I_{k,\nu}^{o,\pm}(z;w).
\end{multline*}
Using that 
$$
\frac{ 1}{\left( q^2;q^2 \right)_{\infty}}
= \omega_{2h,k} (2z)^{\frac{1}{2}} 
\frac{e^{ \frac{\pi}{12k} \left( \frac{ 1}{2z}  - 2 z \right)}}{\left( q_1^{1/2};q_1^{1/2}  \right)_{\infty}}
$$
now gives the claim of the theorem.
\end{proof}
We next consider the case that $k$ is even. 
Define 
\begin{eqnarray*}
l^o(z;w):=  
\frac{1}{(q^2;q^2)_{ \infty }} 
\sum_{m \in \Z } 
\frac{(-1)^m\, q^{3m^2 +5 m}}{1-e^{2 \pi iw } q^{ 2m+1}},
\end{eqnarray*}
and 
\begin{eqnarray*}
J_{ k,\nu}^{o,\pm} (z;w) := \int_{\R} e^{-\frac{3 \pi z}{2k} (x^2-2x)}  
H_w^o \left(  
\frac{2 \pi i}{k} (2 \nu+1)  - \frac{2 \pi zx}{k}  \mp \frac{2 \pi i }{3k} 
\right) dx.
\end{eqnarray*}
\begin{theorem} \label{transtheoremodd2}
Assume that  $k$ is even and $h h' \equiv - 1 \pmod {4k}$. Then we have 
\begin{multline*} 
R^o(z;w) =
i^{-h'} z^{-\frac{1}{2}}   
e^{ - \frac{2\pi z}{3k} - \frac{10 \pi}{3kz}+\frac{\pi i}{2k}(h+7h')}
\, e^{2 \pi i w + \frac{3 \pi  w^2 k}{2z} - \frac{2 \pi w}{z}} \omega_{h,\frac{k}{2}}\, l^o\left(q_1;\frac{iw}{z} \right) \\
+ \frac{1}{k} (-1)^{ \frac{1}{2} (h'+1)}\, e^{ - \pi i w}  \,
\omega_{h,\frac{k}{2}} \, z^{\frac{1}{2}} \,
e^{ \frac{\pi i}{2k} (h-3h'-6) - \frac{2\pi z}{3k} }
 \\
\sum_{ \substack{ \nu \pmod { k/2}  \\ \pm }} 
(-1)^{ \nu} e^{\frac{\pi i}{2k} ( h'( - 12 \nu - 12 \nu^2 \pm4 \nu \pm2 ) -12 \nu \pm 2 )}\, 
J_{ k,\nu}^{o,\pm}(z;w).
\end{multline*}
\end{theorem}
\begin{proof}
Via analytic continuation, we may assume that $w \in \C$ with  Re$(w) \not=0$  sufficiently small and $z>0$ is real.
We use the notation from above and write in   (\ref{tildefunction}) $n = \nu + \frac{k}{2}  m$,
 with $\nu$ running modulo $\frac{k}{2}$ and 
  $m \in \Z$.
This gives that   
\begin{multline} \label{poisson}
\widetilde{R}^o(z;w)=
e^{- \pi i w} \, q^{- \frac{1}{2}} 
\sum_{\nu \pmod{ k/2}} 
(-1)^{\nu} \, 
e^{\frac{6 \pi i h \nu^2}{k} }
\sum_{m \in \Z} (-1)^m 
e^{-\frac{6 \pi z}{k} \left(\nu+ \frac{km}{2}\right)^2}\\
\, H_w^o \left(   \frac{2 \pi i h}{k}(2 \nu+1) - \frac{2 \pi z}{k}  \left(  2\nu +km+1  \right)     \right).
\end{multline} 
We   apply 
Poisson summation  on the inner sum and get,  changing $x$ into $2 \nu + kx$,
 \begin{multline*}
 \widetilde{R}^o(z;w) =
e^{- \pi i w} \, q^{- \frac{1}{2}} \frac{1}{k} 
\sum_{\nu \pmod{ k/2}} 
(-1)^{\nu} \, 
e^{\frac{6 \pi i h \nu^2}{k} }
\sum_{n \in \Z} 
\int_{\R} 
e^{-\frac{3 \pi z x^2}{2k} }\\
\, H_w^o \left(   \frac{2 \pi i h}{k}(2 \nu+1) - \frac{2 \pi z}{k}  ( x+1  )     \right)
e^{ \frac{\pi i}{k }  ( 2n+1) (x- 2 \nu) } dx
.
\end{multline*} 
In the part of the sum on $n$ with $n<0$, we make the change of variables $n \mapsto - n-1$, $x \mapsto -x$, $\nu \mapsto - \nu$, and then $x \mapsto x \pm 1 $ to get 
 \begin{multline*} 
 \widetilde{R}^o(z;w)  =
e^{- \pi i w - \frac{3 \pi z}{2k} } \, q^{- \frac{1}{2}} \frac{1}{k} 
\sum_{\nu \pmod{ k/2}} 
(-1)^{\nu} \, 
e^{\frac{6 \pi i h \nu^2}{k} }
\sum_{\substack{n \geq 0 \\ \pm }} 
\int_{\R} 
e^{-\frac{3 \pi z ( x^2 \mp 2x  )}{2k} }\\
\, H_w^o \left(   \frac{2 \pi i h}{k}(\pm 2 \nu+1) \mp  \frac{2 \pi z x}{k}      \right)
e^{ \frac{\pi i}{k }  ( 2n+1) (x \mp1 - 2 \nu) }\, dx.
\end{multline*}   
We next shift the path of  integration though 
$$
\omega_{n} := \frac{(2n+1)}{3z}.
$$
Using the function $S_w^{\pm}(x)$ defined in (\ref{Sw}), we see that poles of the integrand can only lie in points
\begin{equation} \label{xm3}
x_m^{\pm } := \frac{i}{z} \left(m \pm wk\right).
\end{equation}
Moreover we see that a non-trivial 
residue can only occur if $m$ is odd, and  for fixed  odd $m$, there is at most for one $\nu$ modulo $k/2$ a nontrivial residue.  This can be chosen as 
\begin{equation} \label{nm3}
\nu_m^{\pm} := \frac{1}{2} ( - h' m \mp 1).
\end{equation}
Thus poles may at most occur for $x_m$ and $\nu_m$ as in (\ref{xm3}) and (\ref{nm3}).
Shifting  the path of integration  through 
$\omega_n$, we have to take those $x_m^{\pm}$ into account with $m>0$, $m$ odd,  and 
$n \geq \frac{3m\pm1}{2}$. 
Denoting by $\lambda_{n,m}^{\pm}$ the residue of each  summand, we get
\begin{eqnarray*}
 \widetilde{R}^o(z;w)  =
\sum_1 + \sum_2,
\end{eqnarray*} 
where 
\begin{eqnarray*}
\sum_1 =  \frac{2 \pi i }{k}\, e^{ - \pi i w -  \frac{3 \pi z}{2k}}\, q^{- \frac{1}{2} }
\sum_{ \substack{ m \geq 1, \text{ odd} \\ \pm } }\quad 
\sum_{ n \geq \frac{1}{2} (3m \pm 1)}\, \lambda_{n,m}^{ \pm},
\end{eqnarray*}
 \begin{multline*} 
 \sum_2  =
e^{- \pi i w - \frac{3 \pi z}{2k} } \, q^{- \frac{1}{2}} \frac{1}{k} 
\sum_{\nu \pmod{ k/2}} 
(-1)^{\nu} \, 
e^{\frac{6 \pi i h \nu^2}{k} }
\sum_{\substack{n \in \N_0 \\ \pm }} 
\int_{\R+\omega_n} 
e^{-\frac{3 \pi z}{2k}( x^2 \mp 2x  ) }\\
\, H_w^o \left(   \frac{2 \pi i h}{k}(\pm 2 \nu+1) \mp  \frac{2 \pi z x}{k}      \right)
e^{ \frac{\pi i}{k }  ( 2n+1) (x \mp1 - 2 \nu) }\, dx.
\end{multline*}   
We first consider $\sum_1$. We have
\begin{eqnarray*}
\lambda_{n,m}^{\pm} 
=\pm  \frac{(-1)^{ \nu_m^{ \pm}} k  e^{\frac{6 \pi i h \left( \nu_m^{\pm } \right)^2}{k} - \frac{3 \pi z}{2k} \left(x_m^{ \pm^2}  \mp  2  x_m^{\pm }\right) + \frac{\pi i}{k} (2n+1)( x_m^{\pm}  \mp  1 - 2 \nu_m^{ \pm}   ) }}{2 \pi z  \cosh \left( \mp \pi i w -\frac{\pi i h}{k} (2 \nu_m^{\pm}  \pm1 ) + \frac{\pi zx_m^{\pm}}{k}\right)}.
\end{eqnarray*}
From this one  sees  that
$$
\lambda_{n+1,m}^{\pm}
 = e^{\frac{2 \pi i }{k}  (x_m^{\pm} \mp 1 - 2 \nu_m^{\pm}) }\, \lambda_{n,m}^{\pm}.
$$
A lengthy calculation  gives
\begin{eqnarray*}
\sum_1 =  \frac{1}{z}  e^{2 \pi i w+ \frac{3 \pi w^2 k}{2z}}\, q_1^{\frac{1}{4} }
\left(
\sum_{ \substack{m\geq 1 \text{ odd} \\ \pm }}
(-1)^{ - \frac{h'm}{2}}
\frac{q_1^{\frac{3}{4}  m^2 + \frac{1}{2} ( 1 \pm 1) m }\, e^{-\frac{\pi w}{z}(1\pm1)}}{1-q_1^{m} \, e^{\mp \frac{ 2 \pi  w}{z}}}
\right).
\end{eqnarray*}
In the sum for the $-$-sign  we change $m \mapsto -m$  to obtain 
  a sum over all $m \in \Z$.
Then changing $m \mapsto 2m+1$ gives
\begin{eqnarray*}
\sum_1 =     \frac{i^{ -h'}}{z}  e^{ \frac{3 \pi w^2 k}{2z} + 2 \pi i w - \frac{2 \pi w}{z}  } q^{\frac{1}{4}}  
q_1^{ \frac{7}{4} } \, (q_1^2;q_1^2)_{\infty} \cdot l^o \left(q_1; \frac{iw}{z} \right).
\end{eqnarray*}
We next consider $\sum_2$. In the same way as before we see that
\begin{multline*}
\sum_2= \frac{2}{k} \, e^{ - \pi i w} \, q^{ - \frac{1}{2}}  \sum_{\nu \pmod{k/2} } 
(-1)^{\nu}\, e^{ \frac{6 \pi i h \nu^2}{k}} \\
\sum_{n \in \Z}  
\int_{ \R + \frac{\omega_n}{2}} 
e^{- \frac{6 \pi zx^2}{k}}
H_w^o \left( \frac{2 \pi i h}{k} (  2 \nu+1) -  \frac{2\pi z(2x+1)}{k}   \right)\, 
e^{ \frac{2 \pi i}{k}(2n+1) ( x  -  \nu)} \,dx.
\end{multline*}
Substituting $x \mapsto x + \frac{\omega_n}{2}$ and writing $n = 3p+ \delta$
 with $p \in \Z$  and $\delta \in \{ 0,\pm1\}$ gives
 \begin{multline*}
\sum_2= \frac{2}{k} \, 
e^{ - \pi i w} \, q^{ - \frac{1}{2}}  \sum_{\nu \pmod{k/2} } 
(-1)^{\nu}\, e^{  \frac{6 \pi i h \nu^2}{k}} 
\sum_{p,\delta}  
e^{ - \frac{\pi \left(6p+2 \delta+1\right)^2}{6kz}- \frac{2 \pi i}{k}( 6p+ 2 \delta +1)\nu}
 \\
\int_{ \R} 
e^{- \frac{6 \pi zx^2}{k}} 
H_w^o \left( 
\frac{2 \pi i }{k}
 \left(  h \left(2 \nu+1 \right) - 2 p \right)-  \frac{2\pi z(2x+1)}{k}  - \frac{2 \pi i }{3k} \left( 2 \delta+1\right) 
\right)\, 
 \,dx.
\end{multline*}
Changing $\nu \mapsto \frac{1}{2} ( -  h'  ( 2 \nu + 2 p +1) -1 )$  and $x \mapsto 2 x+1$ gives after a lengthy calculation 
\begin{multline*}
\sum_2 
= e^{- \pi i w } q^{ \frac{1}{4}} \frac{1}{k} (-1)^{ \frac{1}{2}(h'+1)}  e^{-\frac{3 \pi i h'}{2k} -\frac{3\pi i}{2k}} 
\sum_{\substack{\nu \pmod{ k/2}  \\ \pm}}
(-1)^{\nu} 
\\
e^{\frac{\pi i}{2k} h'\left(- 12 \nu - 12 \nu^2 \pm 4 \nu \pm 2    \right)}
e^{\frac{\pi i}{2k}(- 12 \nu \pm2 )} 
\sum_{p} (-1)^p \, q_1^{\frac{1}{12}(6p+2 \delta+1)^2 }
 J_{k,\nu}^{o,\pm}(z;w).
\end{multline*} 
The sum over $p$   vanishes for $\delta=1$. If $\delta=0,-1$, then it equals 
$\eta\left( \frac{4 \pi i}{k}(h'+ \frac{i}{z} )\right)$.
Using 
\begin{eqnarray*}
\frac{1}{(q^2;q^2)_{\infty}} 
= \frac{\omega_{h,\frac{k}{2}} \, z^{\frac{1}{2}} \, e^{\frac{\pi}{6k} (z^{-1}-z) }}{(q_1^2;q_1^2)_{\infty}}
\end{eqnarray*}
now easily gives the claim.
\end{proof}
We next conclude  from Theorems  \ref{transtheoremodd1}  and 
 \ref{transtheoremodd2}
 a transformation law for $R_2^o(q)$. For this  
 let 
 \begin{eqnarray*}
g^o(q)&:=& \frac{1}{(q;q)_{\infty} } \sum_{m \in \Z} (-1)^m \frac{q^{ \frac{1}{2}(3m^2+m)}}{(1+q^m)},\\
h^o(q)&:=&  \frac{1}{(q;q)_{\infty} } 
\sum_{m \in \Z} (-1)^m
\frac{q^{\frac{1}{2} (3m^2+m) } (1-3q^m)}{(1+q^m)^3},\\
I_{k,\nu}^{o,\pm}(z) &:=& \int_{\R}
\frac{  e^{ - \frac{3 \pi z x^2}{2k} + \frac{\pi zx}{k}}  }{
\sinh^3 \left(\mp \frac{\pi i \nu}{k} + \frac{\pi zx}{k} \pm \frac{\pi i}{6k} \right)
}\, dx,\\
m^o(q)&:=&  \frac{1}{(q^2;q^2)_{\infty} }   \sum_{m\in \Z} (-1)^m \frac{q^{3m^2+5m}}{1-q^{2m+1}},\\
J_{k,\nu}^{o,\pm}(z) &:=& \int_{\R}
\frac{e^{ - \frac{3 \pi z x^2}{2k} + \frac{ \pi z x }{k}}}{\sinh^3 \left(- \frac{\pi i}{k} (2 \nu+1) + \frac{\pi z x}{k} \pm \frac{\pi i }{3k} \right)}\, dx.
\end{eqnarray*}
\begin{corollary} \label{corollarytransfodd}
\begin{enumerate}
\item 
If $k$ is odd, then 
\begin{multline*}
R_2^o(q) = - \frac{1}{\sqrt{2}}
 \omega_{2h,k} \,  \
 e^{\frac{\pi}{24kz} - \frac{2 \pi z}{3k}+ \frac{\pi i h}{2k}  + \frac{\pi i}{2}  (k+1+3hk)}
\left( -
 \frac{3k}{8\,\pi\, z^{\frac{3}{2}}} g^o\left(  q_1^{\frac{1}{2} }\right) - \frac{1}{16z^{\frac{5}{2}} }
 h^o\left(q_1^{\frac{1}{2}} \right)
\right) \\ 
 - \frac{1}{8 \sqrt{2}k}\, e^{\frac{\pi i}{2} k(3h+1) + \frac{\pi i h}{2k}- \frac{2 \pi z}{3k}}\, \omega_{2h,k}\, z^{\frac{1}{2} }\,  
 \sum_{\substack{ \nu \pmod k\\ \pm} } \pm e^{\frac{\pi i}{2k}  \left(  h' \left(- 3 \nu^2 + \nu \right) \pm \frac{1}{3}  \mp 2 \nu\right)}\, I_{k,\nu}^{o,\pm}(z)
 .
\end{multline*}
\item
If $k$ is even, then 
\begin{multline*}
R^o_2(q) = i^{-h'}   e^{ - \frac{2\pi z}{3k} - \frac{10 \pi}{3kz} + \frac{\pi i}{2k}(h+7h') }\,  
\omega_{h,\frac{k}{2}} 
\left( 
- \frac{3k}{8\pi z^{\frac{3}{2}}} m^o(q_1)  - \frac{1}{z^{\frac{5}{2}}} q_1^{-2} R_2^o(q_1)
\right)\\
+ \frac{1}{8k} (-1)^{ \frac{1}{2}(1+h')}\, \omega_{h,\frac{k}{2}}\, z^{\frac{1}{2}}
  e^{\frac{\pi i}{2k} (h-3h') -\frac{2 \pi z}{3k}  -\frac{\pi i}{k}  }\\
\sum_{\substack{\nu \pmod{k/2}\\ \pm }}
(-1)^{\nu}\,e^{ \frac{\pi i h'}{k} \left(- 6 \nu - 6 \nu^2 \pm 2 \nu \pm 1    \right)- \frac{2 \pi i \nu}{k} 
 \pm \frac{\pi i}{3k}} \, J_{k,\nu}^{o,\pm}(z).
\end{multline*}
\end{enumerate}
\end{corollary}
\begin{proof}
We have 
\begin{eqnarray*}
L^o
\left(   \frac{
e^{
\frac{3 \pi kw^2 }{2z} + 2 \pi iw - \frac{\pi w}{2z}}
}{1+ e^{ - \frac{\pi w}{z} }\, q^m}
\right)
= \frac{i(1-q^m)}{8z(1+q^m)^2}  - \frac{3k}{8\pi z(1+q^m)}
- \frac{1-6q^m+q^{2m}}{32z^2(1+q^m)^3}.
\end{eqnarray*}
Moreover
\begin{eqnarray*}
\sum_{m \in \Z} (-1)^m \frac{q^{ \frac{1}{2}(3m^2+m)} (1-q^m)}{(1+q^m)^2} =0,
\end{eqnarray*}
since the $m$th and $-m$th term cancel, and 
\begin{eqnarray*}
\sum_{m \in \Z}(-1)^m \frac{q^{ \frac{1}{2} (3m^2+m)}\left(1- 6q^m + q^{2m} \right)}{(1+q^m)^3}
= 2 \sum_{m \in \Z} (-1)^m \frac{q^{\frac{1}{2}(3m^2+m) } (1-3q^m)}{(1+q^m)^3}.
\end{eqnarray*}
Furthermore 
\begin{eqnarray} \label{diff}
L^o \left(e^{- \pi i w} \, H_w^o(x)\right) =   \frac{e^x}{8 \sinh^3\left( - \frac{x}{2}\right)}.
\end{eqnarray}
Now (i) follows easily using Theorem \ref{transtheoremodd1}.

To see (ii), we observe that 
\begin{eqnarray*}
L^o  \left(   \frac{
e^{
\frac{3 \pi kw^2 }{2z} + 2 \pi iw - \frac{2\pi w}{z}}
}{1 - e^{ - \frac{2 \pi w}{z} }\, q^{2m+1}}
\right) 
&=& 
\frac{i}{2z(1-q^{2m+1 })^2} 
- \frac{3k}{8 \pi z (1-q^{2m+1} )} 
- \frac{1+q^{2m+1}}{2z^2(1-q^{2m+1})^3},
\end{eqnarray*}
use again (\ref{diff}), and  
\begin{eqnarray*}
\sum_{m \in \Z} \frac{(-1)^m\, q^{3m^2 + 5m}}{(1-q^{2m+1})^2} =0,
\end{eqnarray*}
since the $m$th and $-(m+1)$th  term cancel.
Then (ii) follows from Theorem  \ref{transtheoremodd2}.
\end{proof}
\section{Proof of Theorem \ref{maintheorem3}}  \label{AsympoddSection}
Here we show asymptotics for $\eta_2^o(n)$. As in Section \ref{TRANS}, we   estimate the Mordell type integrals  occurring in the transformation  law of $R^o(q)$  and Kloosterman sums. 
We assume the notation from above.
\begin{lemma} \label{lemmaintest2}
\begin{enumerate}
\item
If $k$ is odd, then 
\begin{eqnarray*}
z^{\frac{1}{2} } I_{k,\nu}^{o,\pm}(z)
 \ll \frac{n^{ \frac{1}{4}}}{\left\{ \frac{\nu}{k}- \frac{1}{6k}  \right\}^3} .
\end{eqnarray*}
\item
If $k$ is even, then 
\begin{eqnarray*}
z^{\frac{1}{2} } 
J_{k,\nu}^{o,\pm}(z)
 \ll \frac{n^{ \frac{1}{4}}}{\left\{ \frac{2\nu+1}{k}\mp \frac{1}{3k}  \right\}^3} .
\end{eqnarray*}\end{enumerate}
\end{lemma}
\begin{lemma} \label{lemmakloostest2}
Let $l,n\in \Z$.
\begin{enumerate}
\item 
If $k$ is odd, and $hh' \equiv - 1 \pmod k$ with $h'$ even, then 
\begin{eqnarray} \label{Kloostodd1}
\sum_{ \substack{h \pmod{k}^* \\ 
\sigma_1 \leq Dh' \leq \sigma_2 }}
\omega_{2h,k} \, e^{ \frac{3 \pi i hk}{2}}
e^{\frac{\pi i}{2k} ( h(1+4n) + 2 lh' ) }
\ll (3n+1,k)^{\frac{1}{2} + \epsilon}\, k^{\frac{1}{2} + \epsilon }.
\end{eqnarray}
\item
If $k$  is even and $hh' \equiv -1 \pmod {4k}$, then
\begin{eqnarray} \label{Kloostodd2}
\sum_{ \substack{h \pmod{k}^* \\ \sigma_1 \leq Dh' \leq \sigma_2 }} 
i^{-h'} \,  
 \omega_{h, \frac{k}{2}}
e^{\frac{\pi i}{2k} ((1+4n)h+(7+4l)h' )} \ll
(96n+25,k)^{\frac{1}{2}}\, k^{\frac{1}{2} + \epsilon  },
\end{eqnarray}
\begin{eqnarray} \label{Kloostodd3}
\sum_{ \substack{h \pmod{k}^* \\ \sigma_1 \leq Dh' \leq \sigma_2 }} 
i^{h'} \,  
 \omega_{h, \frac{k}{2}}
e^{\frac{\pi i}{2k} ((1+4n)h+(-1+4l)h' )} \ll
(96n+25,k)^{\frac{1}{2}}\, k^{\frac{1}{2} + \epsilon  }.
\end{eqnarray}
\end{enumerate}
\end{lemma}
\begin{proof}
(i) 
Using that $h'$  is even, one can see  that  the left-hand side of (\ref{Kloostodd1}) is  well-defined.
We only show  the estimate for the full modulus $k$, the restriction to   $\sigma_1 \leq Dh' \leq \sigma_2 $ can be concluded as in \cite{Le}. 
We change $h$ into $\bar 2 h$ and $h'$ into $2h'$.
This leads to the sum
\begin{eqnarray} \label{Kloostodd4}
\sum_{h \pmod{k}^*} \omega_{h,k} \, e^{\frac{2\pi i}{k}  \left( \left( \bar 2 n + \frac{1}{4} (3k^2+1) \bar{2} \right) h + l h'  \right)  }.
\end{eqnarray}
By \cite{An1}, (\ref{Kloostodd4}) can be estimated against   
 \begin{eqnarray*}
  \left( 24\left(  \bar 2  n + \frac{1}{4} \left(3k^2+1\right) \bar{2} \right)+1,k \right)^{\frac{1}{2} } 
  k^{\frac{1}{2} + \epsilon }
 = (3n+1,k)^{\frac{1}{2} }\, k^{\frac{1}{2} + \epsilon  }.
 \end{eqnarray*}
 (ii)  We only show (\ref{Kloostodd2}), since (\ref{Kloostodd3}) can be proven similarly.
 We first observe that (\ref{Kloostodd2})
 is well-defined since changing $h \mapsto h+k$ implies that $h' \mapsto h' + rk$ with $r \equiv k+1 \pmod 4$.
 We then change the sum into a sum $\pmod{4k}$. 
 This gives  
  \begin{eqnarray*}
 \frac{1}{4} \sum_{ h \pmod {4k}^*}
 \omega_{h,\frac{k}{2} }\, e^{\frac{2\pi i}{4k} ( h(4n+1) + (4l-7-k) h'     )} 
 \ll (96n+25,k)^{ \frac{1}{2} } k^{\frac{1}{2} + \epsilon }
  \end{eqnarray*}
  by \cite{An1} as claimed .
\end{proof}
\begin{proof}[Proof of Theorem \ref{maintheorem3}]
We  use  the circle method and proceed as in the proof of Theorem \ref{maintheorem1}. 
This gives  
\begin{eqnarray*}
\eta_2^o \left(n \right) = 
\sum_{h,k} 
e^{- \frac{2 \pi i hn}{k}}
\int_{-\vartheta_{h,k}'}^{\vartheta_{h,k}''}
R_2^o \left(e^{\frac{2 \pi i }{k}(h+iz) } \right)  \cdot 
e^{\frac{2 \pi n z}{k}} \ d\Phi.
\end{eqnarray*}
Corollary \ref{corollarytransfodd} gives
\begin{multline*}
\eta_2^o(n) = - 
\frac{i}{\sqrt{2}} \sum_{\substack{k \text{ odd} \\ h}} 
e^{- \frac{2 \pi i hn}{k} + \frac{\pi i  h}{2k}  + \frac{\pi i}{2} (k+3hk)} 
\omega_{2h,k} 
 \int_{ -\vartheta_{h,k }' }^{\vartheta_{h,k}^{''} }
e^{\frac{2 \pi  z}{ k}\left( n- \frac{1}{3} \right) + \frac{\pi}{24 kz}} 
\left( - \frac{3\pi k}{8z^{\frac{3}{2} }}  g^o \left( q_1^{\frac{1}{2} }  \right) 
\right.
\\
\left.
- \frac{1}{16z^{\frac{5}{2} } }  h^o\left( q_1^{ \frac{1}{2}} \right)\right) d \Phi
+ \sum_{\substack{k \text{ even} \\ h}} i^{-h' } e^{-\frac{2 \pi i hn}{k} +\frac{\pi i}{2k}(h+7h')} \omega_{h, \frac{k}{2}} 
 \int_{ -\vartheta_{h,k }' }^{\vartheta_{h,k}^{''} }e^{\frac{2 \pi  z}{ k}\left( n- \frac{1}{3} \right) - \frac{10\pi}{3kz}}
\left(-\frac{3 k}{8 \pi z^{\frac{3}{2} }} m^o(q_1)
\right. \\
\left.
 - \frac{1}{ z^{\frac{5}{2} } } q_1^{-2}  R_2^o(q_1) \right)
- \frac{1}{8 \sqrt{2}} 
\sum_{ \substack{k \text{ odd}  \\ h } } 
\frac{\omega_{2h,k }}{k}\, e^{-\frac{2 \pi i hn}{k} }\,e^{\frac{\pi i }{2} k(3h+1)  + \frac{\pi i h}{2k} }  
\sum_{\substack{ \nu \pmod {k} \\ \pm}  } \pm
e^{ \frac{\pi i}{2k} \left( h' \left( -3 \nu^2 +\nu \right) \pm \frac{1}{3} \mp 2 \nu  \right)} \\
 \int_{ -\vartheta_{h,k }' }^{\vartheta_{h,k}^{''} }
 z^{\frac{1}{2} }
e^{ \frac{2 \pi z}{k}\left( n- \frac{1}{3} \right)} I_{k,\nu}^{o,\pm}(z) d\Phi 
+ \frac{1}{8k} \sum_{ \substack{k \text{ even}  \\ h } } 
(-1)^{\frac{1}{2} (1+h') } \omega_{h,\frac{k}{2}} \,
e^{\frac{\pi i}{2k} (h-3h') -\frac{\pi i}{k}} \, e^{-\frac{2 \pi  i hn}{k} }
\sum_{\substack{ \nu \pmod {k/2} \\ \pm}  } 
(-1)^{\nu} \\
e^{\frac{\pi i}{k} (h' (- 6 \nu - 6 \nu^2 \pm 2 \nu \pm 1) ) - \frac{2\pi i \nu}{k}\pm \frac{\pi i }{3k} }
 \int_{ -\vartheta_{h,k }' }^{\vartheta_{h,k}^{''} }
 z^{\frac{1}{2} }
e^{ \frac{2 \pi z}{k}\left( n- \frac{1}{3} \right)} J_{k,\nu}^{o,\pm}(z) \, d\Phi
.
\end{multline*}
We denote the occurring sums by 
$\sum_1,\sum_2,\sum_3$, and $\sum_4$. 
As in the proof of Theorem \ref{maintheorem1} we obtain,
 using Lemma \ref{lemmakloostest2},
\begin{eqnarray*}
\sum_1 + \sum_2 
=- \frac{i}{\sqrt{2}} \sum_{k \text{ odd} }
A_k^o(n)    
 \int_{ -\frac{1}{kN}}^{\frac{1}{kN} }
e^{\frac{2 \pi  z}{ k}\left( n- \frac{1}{3} \right) + \frac{\pi}{24 kz}} 
\left( - \frac{3\pi k}{16z^{\frac{3}{2} }}   + \frac{1}{64z^{\frac{5}{2} } }   \right) d \Phi
+ O(n^{2 + \epsilon}).
\end{eqnarray*}
We consider integrals of the form
\begin{eqnarray*}
I_{k,r} :=
\int_{-\frac{1}{kN} }^{ \frac{1}{kN}}
z^r e^{\frac{2 \pi}{k} \left( z \left( n-\frac{1}{3}\right)+ \frac{1}{48z} \right) } d\Phi
\end{eqnarray*}
with $r \in \{-\frac{5}{2},-\frac{3}{2} \}$.
As in the proof of Theorem \ref{maintheorem1} this equals
\begin{eqnarray*}
I_{k,r} = \frac{1}{ki}
\int_{\Gamma} z^r   e^{\frac{2 \pi}{k} \left( z \left( n-\frac{1}{3}\right)+ \frac{1}{48z} \right) } 
dz + O(n^{\frac{3}{4}}).
\end{eqnarray*}
Making the substitution $t = \frac{\pi}{24 kz}$ gives
\begin{eqnarray*}
I_{k,r}=
\frac{2\pi}{k} \left( \frac{\pi}{24 k}\right)^{1+r}
\frac{1}{2 \pi i}\int_{\gamma - i \infty }^{\gamma+ i \infty}
t^{-2 -r } e^{ t+ \frac{\alpha}{t}} dt 
+ O \left( n^{ \frac{3}{4}}\right),
\end{eqnarray*}
where $\gamma \in \R$ and $\alpha =  \frac{\pi^2}{36k^2} (3n-1) $.
Now the integrals can be computed as before and lead to Bessel functions of order $-\frac{1}{2}$ and $-\frac{3}{2}$.
More precisely, we obtain 
that  $\sum_1 +\sum_2$ equals
\begin{eqnarray*}
 -\frac{i}{\sqrt{2}} \sum_{k \text{ odd} } 
A_k^o(n) 
\left(
-\frac{3\pi^2(3n-1)^{\frac{1}{4}}}{4}  I_{-\frac{1}{2} } \left( \frac{\pi}{3k} \sqrt{3n-1}\right)
+ \frac{\pi (3n-1)^{\frac{3}{4}}}{4k}  I_{-\frac{3}{2} } \left( \frac{\pi}{3k} \sqrt{3n-1}\right)
\right) + O(n^{2+\epsilon}).
\end{eqnarray*}
$\sum_3$ and $\sum_4$ can be estimated against $O\left(n^{ 2+\epsilon} \right)$ as in the proof of Theorem \ref{maintheorem1} using Lemma \ref{lemmakloostest2} and Lemma \ref{lemmaintest2}.
\end{proof}

\end{document}